\documentclass[11pt]{article}

\usepackage{latexsym,amssymb}
\usepackage{amsmath}
\usepackage{mathrsfs}
\usepackage{bussproofs}
\usepackage{tikz}

\usepackage{enumerate}
\usepackage{multicol}
\usepackage{color}

\newtheorem{defi}{Definition}[section]
\newtheorem{theo}[defi]{Theorem}
\newtheorem{lem}[defi]{Lemma}
\newtheorem{prop}[defi]{Proposition}
\newtheorem{rem}[defi]{Remark}

\newtheorem{cor}[defi]{Corollary}

\newcommand{\negr}[1]{\boldsymbol{#1}}

\newenvironment{dem}{\noindent\bf Proof. \rm}{\hfill $\negr{\blacksquare}$}

\def\Hom{\mathop{\rm Hom}\nolimits}
\def\T{\mathop{\rm T}\nolimits}
\def\F{\mathop{\rm F}\nolimits}

\def\n{\mathbf{n}}
\def\b{\mathbf{b}}
\def\0{\mathbf{0}}
\def\1{\mathbf{1}}

\newcommand{\tml}{{\cal TML}}

\title{Normal proofs and tableaux  for the Font-Rius tetravalent modal logic}

\author{
    Marcelo E. Coniglio $^\textup{\scriptsize a}$
    \and
    Mart\'{\i}n Figallo
  $^\textup{\scriptsize b}$
}

\date{
    $^\textup{\scriptsize a}$\textit{\small CLE and Department of Philosophy. State University of Campinas. Campinas, Brazil}
    \\
    $^\textup{\scriptsize b}$\textit{\small Departamento de Matem\'atica and Instituto de Matem\'atica (INMABB). Universidad Nacional del Sur. Bah\'{\i}a Blanca, Argentina}
}

\begin{document}

\maketitle

\begin{abstract}
{\em Tetravalent modal logic} ($\cal TML$) was introduced by Font and Rius in 2000; and it is an expansion of the Belnap-Dunn four--valued logic $\cal FOUR$, a logical system that is well--known for the many applications it has been found in several fields.  Besides, $\cal TML$ is the logic that preserve degrees of truth with respect to Monteiro's tetravalent modal algebras. Among other things, Font and Rius showed that $\cal TML$ has a strongly adequate sequent system, but unfortunately this system does not enjoy the cut-elimination property. However, in  \cite{MF}, it was presented a sequent system for $\cal TML$ with the cut--elimination property. Besides, in this same work, it was also presented a sound and complete natural deduction system for this logic.

In the present work, we continue with the study of $\tml$ under a proof-theoretic perspective. In first place, we show that the natural deduction system given in \cite{MF} admits a normalization theorem. In second place, taking advantage of the contrapositive implication for the tetravalent modal algebras introduced in \cite{AVF3}, we define a decidable tableau system adequate to check validity in the logic $\tml$.  Finally, we provide a sound and complete tableau system for $\tml$ in the original language.These two tableau systems constitute new (proof-theoretic) decision procedures for checking validity in the variety of tetravalent modal algebras.
\end{abstract}

\

\section{Introduction}\label{s1}
The class {\bf TMA} of tetravalent modal algebras was first
considered by Antonio Monteiro (1978), and mainly studied by I.
Loureiro, A.V. Figallo, A. Ziliani and P. Landini. 
From Monteiro's point of view,  in the future these algebras would 
give rise to a four-valued modal logic with significant applications in Computer 
Science (see \cite{FR2}). Later on, J.M.
Font and M. Rius were indeed interested in the logics arising from the
algebraic and lattice--theoretical aspects of these algebras.

Although such applications have not yet been developed, the two logics considered in \cite{FR2} are modal expansions
of Belnap-Dunn's four-valued logic, a logical system that is well--known for the many applications it has found in several fields. 
In these logics, the four non-classical epistemic values emerge: $\1$ (true and not false), $\0$ (false and not true), $\n$ (neither true nor false) and $\b$ (both true and false). 
We may think of them as the four possible ways in which an atomic sentence $P$ can
belong to the {\em present state of information} : we were told that (1)  $P$ is true (and were not told that $P$ is false); (2) $P$
is false (and were not told that $P$ is true); (3) $P$ is both true and false (perhaps from different sources, or in different instants
of time); (4) we were not told anything about the truth value of $P$.
In this interpretation, it makes sense to consider a modal-like unary operator $\square$ of epistemic character, such that for any sentence $P$, the sentence $\square P$ would mean  ``the available information confirms that $P$ is true".
 It is clear that in this setting the sentence $\square P$ can only be true in the case where we have some information saying that $P$  is true and we have no information saying that $P$ is false, while it is simply false in all other cases (i.e., lack of information or at least some information saying that $P$ is false, disregarding whether at the same time some other information says that $P$ is true); that is, on the set $\{\0, \n, \b, \1\}$ of epistemic values this operator must be defined as $\square \1 = \1$ and $\square \n = \square \b =\square \0=\0$ . This is exactly the algebra that generates the variety of tetravalent modal algebras (TMAs). 

In \cite{FR2}, Font and Rius studied two logics related to TMAs. One of them is obtained by following the usual ``preserving truth''  scheme,  taking  $\{\1\}$ as designated  set, that is, $\psi$ follows from $\psi_1, \dots, \psi_n$ in  this  logic  when every interpretation  that  sends  all  the  $\psi_i$ to  $\1$  also sends  $\psi$ to  $1$. The  other logic, denoted by ${\cal TML}$ (the logic we are interested in), is defined by  using  the  {\em preserving  degrees of truth}  scheme, that is, $\psi$  follows  from  $\psi_1, \dots, \psi_n$  when  every interpretation that assigns to  $\psi$ a value that is greater  or equal than the value it assigns to  the conjunction of  the  $\psi_i$'s.  These authors proved that ${\cal TML}$  is  not  algebraizable  in  the  sense of Blok and Pigozzi, but it is  {\em finitely equivalential}  and  {\em protoalgebraic}. However, they confirm that its algebraic  counterpart  is also the  class of TMAs: but the  connection between the logic and the algebras is not so good  as in the first logic.  As a compensation,  this logic has a better  proof-theoretic  behavior,  since  it  has  a {\em strongly  adequate  Gentzen 
calculus} (Theorems  3.6 and 3.19 of \cite{FR2}). 

In \cite{FR2}, it was proved that $\cal TML$ can be characterized as a matrix logic in terms of two logical matrices, but later, in \cite{FR2}, it was proved that $\cal TML$ can be determined by a single logical matrix (see Proposition \ref{prop2.5} below). Besides, taking profit of the contrapositive implication introduced by A. V. Figallo and P. Landini (\cite{AVF3}), a  sound and complete Hilbert-style calculus for this logic was presented. Finally, the paraconsistent character of  $\cal TML$ was also studied in \cite{CF} from the point of view of the \emph{Logics
of Formal Inconsistency}, introduced by W. Carnielli and J. Marcos in~\cite{Tax} and afterward developed in~\cite{WCMCJM}. 

\

In this work, we continue with the study of $\tml$ from a syntactic point of view. In Section \ref{s2}, we recall the natural deduction system presented in \cite{MF} with introduction and elimination rules for every connective; and we show that every deduction ${\cal D}$ in this system reduces to a normal form. In Section \ref{s3}, we present $\tml$ in a different signature containing the contrapositive implication defined in \cite{AVF3}. In Section \ref{s4}, by adapting the general technique introduced in~\cite{CCCM}, a decidable tableau system is defined, in the language $\{\neg, \succ\}$, proving the corresponding soundness and completeness theorems. Finally, in Section \ref{s6}, we provide a new tableau system for $\tml$ in the original language.

\section{Preliminaries}%\label{s1}

Recall that, a {\em De Morgan algebra} is a structure $\langle
A, \wedge, \vee, \neg, \0 \rangle$ such that $\langle A, \wedge, \vee, \0\rangle$ 
is a bounded distributive lattice and
$\neg$ is a De Morgan negation, i.e., an involution that
additionally satisfies De Morgan's laws: for every $a, b\in A$ \, $\neg \neg a = a$, and $\neg (a \vee b) = \neg a \wedge  \neg b$.\\
A {\em tetravalent modal algebra} (TMA) is an algebra
$\mathbb{A}=\langle A, \wedge, \vee, \neg, \square, \0\rangle$ of
type $(2,2,1,1,0)$ such that its non-modal reduct
$\langle A, \wedge, \vee, \neg, \0\rangle$ is a De
Morgan algebra and the unary operation $\square$ satisfies the identities
all $a\in A$, the two following axioms:
 $$\square a \wedge \neg a \approx \0,$$ 
 $$\neg \square a \wedge a \approx \neg a \wedge a.$$ 
Every TMA $\mathbb{A}$ has a top element $\1$ which is defined as $\neg \0$. 
These algebras were studied mainly by I. Loureiro (\cite{L1,L2}), and also by A. V.  Figallo, P. Landini (\cite{AVF3}) and A. Ziliani (\cite{AVF2}), at the suggestion of the late A. Monteiro (see \cite{FR2}).
The class of all tetravalent modal algebras constitute
a variety which is denoted by {\bf TMA}. Let $M_4=\{\0,\n, \b, \1\}$ and consider the lattice given by the following Hasse diagram 
 
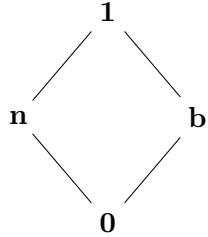
\begin{center}
\begin{tikzpicture}[scale=.7]
%\tikzstyle{every node}=[draw,circle,fill=white,inner sep=2pt]
  \node (one) at (0,4) [] {$\1$};
  \node (b) at (-1.7,2) [] {$\n$};
  \node (a) at (1.7,2) [] {$\b$};
  \node (zero) at (0,0) [] {$\0$};
  \draw (zero) -- (a) -- (one) -- (b) -- (zero);
\end{tikzpicture} \hspace{2cm}
\end{center}
This is a well-known lattice and it is called ${\bf L4}$ (See \cite{AnBel}, pg. 516.) 
Then, {\bf TMA} is generated by the above  four-element   lattice enriched with two unary operators $\neg$ and $\square$  given by $\neg \n= \n$, $\neg \b=\b$, $\neg \0=\1$ and $\neg \1=\0$ and the unary operator $\square$ is defined as: $\square \n= \square \b = \square \0=\0$ and $\square \1=\1$  (see \cite{FR2}). This tetravalent modal algebra, denoted by  $\mathfrak{M}_{4m}$, has two prime filters, namely, $F_{\tiny \mbox{\bf n}} =\{\n,  \1\}$ and $F_{\tiny \mbox{\bf b}} =\{\b, \1\}$. As we said, $\mathfrak{M}_{4m}$ generates the variety ${\bf TMA}$, i.e., an
equation holds in every  TMA  iff it holds in $\mathfrak{M}_{4m}$. 

\begin{lem}(See \cite{FR2})\label{lem2.1}\label{PropTMA} In every TMA $\mathbb{A}$ and for all $a,b\in A$ the following identities hold.
\begin{multicols}{2}
\begin{enumerate}[\rm (i)]
\item $\neg \square a \vee  a\approx\1$,
\item  $\square \square a \approx  \square a$,
\item $\square a \vee  \neg a \approx a \vee  \neg a $,  
\item $\square(a\wedge b)\approx\square a\wedge \square b$, 
\item $\square a \vee \neg\square a\approx\1$, 
\item $\square (a \vee \square b) \approx \square a \vee \square b$,
\item $\square a \wedge \neg \square a \approx\0$,
\item $\square \neg \square a \approx \neg\square a$ 
\item $\square a \wedge  a \approx a$, 
\item $ a\wedge \square \neg a \approx \0$,
\item $\square \1 \approx \1$,  
\item $\square(\square a\wedge \square b)\approx\square a\wedge \square b$,
\item $\square \0 \approx \0$, 
\item$\square(\square a\vee \square b)\approx\square a\vee \square b$.
\end{enumerate}
\end{multicols}
\end{lem}

Let $\mathscr{L}=\{\vee, \wedge, \neg, \square\}$ be a propositional language. 
From now on, we shall denote by $\mathfrak{Fm}$ the absolutely free algebra of
type (2,2,1,1,0) generated by some denumerable set of variables, i.e., $\mathfrak{Fm} =\langle Fm, \wedge,
\vee, \neg, \square, \bot \rangle$. As usual, we refer to formulas
by lowercase Greek letters $\alpha, \beta, \gamma, \dots$ and so
on; and to finite sets of formulas by uppercase Greek
letters $\Gamma, \Delta,$ etc.. Then

\begin{defi}(\cite{FR2, CF}) \label{logTMA} The tetravalent modal logic ${\cal TML}$
defined over $\mathfrak{Fm}$ is the propositional logic
$\langle Fm, \models_{{\cal TML}}\rangle$ given as
follows: for every finite set $\Gamma \cup\{\alpha\} \subseteq
Fm$, $\Gamma \models_{\cal TML} \alpha$  if and only if, for every
$\mathbb{A} \in {\bf TMA}$ and for every $h \in
\Hom(\mathfrak{Fm}, \mathbb{A})$, $\bigwedge \{ h(\gamma) \ : \
\gamma \in \Gamma\} \leq h(\alpha)$. In particular, $\emptyset
\models_{\cal TML} \alpha$ if and only if $h(\alpha)=\1$ for every
$\mathbb{A} \in {\bf TMA}$ and for every $h \in
\Hom(\mathfrak{Fm}, \mathbb{A})$.
\end{defi}

\begin{rem} Observe that, if $h\in \Hom(\mathfrak{Fm}, \mathbb{A})$ for any $\mathbb{A} \in {\bf TMA}$, we have that $h(\bot)=\0$. This follows from the fact that $\bot$ is the $0$-ary operation in $\mathfrak{Fm}$, $\0$ is the $0$-ary operation in  $\mathbb{A}$ and the definition of homomorphism (in the sense of universal algebra). 
\end{rem}

It is worth mentioning that there is a number of works on modal logics which either share the non-modal fragment with ${\cal TML}$ or have
non-modal fragments which are characterized by the same four-element matrix. Clearly, these logics have some relation to ${\cal TML}$. Some examples of such systems are Priest’s KFDE \cite{Pri},
Belnapian modal logics of Odintsov and Wansing \cite{OdWa01, OdWa02} and modal bilattice logic \cite{RiJuJa}.\\[2mm]
The following result was proved in \cite{CF} and will be useful in the sequel.
\begin{lem}\label{Lemhh'} Let $h \in \Hom(\mathfrak{Fm},\mathfrak{M}_{4m})$,
${\cal V}' \subseteq Var$ \, and \, $h' \in \Hom(\mathfrak{Fm},\mathfrak{M}_{4m})$
such that, for all $p \in {\cal V}'$,\\[4mm]
$h'(p)= \left \{ \begin{tabular}{cl}
$h(p)$ & if $h(p)\in \{\0,\1\}$,\\
$\n$ & if $h(p)=\b$,\\
$\b$ & if $h(p)=\n$\\
\end{tabular}\right.$ .  Then\,  $h'(\alpha)= \left \{ \begin{tabular}{cl}
$h(\alpha)$ & if $h(\alpha)\in \{\0,\1\}$,\\
$\n$ & if $h(\alpha)=\b$,\\
$\b$ & if $h(\alpha)=\n$,\\
\end{tabular}\right.$ for all $\alpha \in \mathfrak{Fm}$ whose variables are in ${\cal V}'$.
\end{lem}

\

\

J. M. Font and M. Rius proved in \cite{FR2} that the tetravalent modal logic $\cal TML$ is a matrix logic. In fact, $\tml$ can be determined by the matrix ${\cal M}_{\n}=\langle\mathfrak{M}_{4m},\{\n, \1\}\rangle$ and simultaneously, it can be determined by ${\cal M}_{\b}=\langle\mathfrak{M}_{4m},\{\b, \1\}\rangle$ (both matrices are isomorphic).

\begin{prop}\label{prop2.5} (\cite{FR2}) $\tml$ is the logic determine by the matrix ${\cal M}_{\b}=\langle\mathfrak{M}_{4m},\{\b, \1\}\rangle$.
\end{prop}

In order to characterize $\cal TML$ syntactically, that is, by means of a syntactical deductive system, it
was introduced in~\cite{FR2} the Gentzen--style system $\mathfrak{G}$. The sequent calculus $\mathfrak{G}$  is single--conclusion, that is, it deals with sequents of the form $\Delta \Rightarrow \alpha$ such that $\Delta \cup \{\alpha\}$ is a finite subset of $Fm$. The axioms and rules of $\mathfrak{G}$ are the following:

\noindent {\bf Axioms}
$$ \mbox{(Structural axiom) \, } \displaystyle {\alpha \Rightarrow \alpha} \hspace{2cm} \mbox{(Modal axiom) \,  } {\Rightarrow \alpha \vee \neg \square \alpha}$$

\noindent {\bf Structural rules}

$$ \mbox{(Weakening) \, } \displaystyle \frac{\Delta \Rightarrow \alpha} {\Delta, \beta \Rightarrow \alpha} \hspace{2cm} \mbox{(Cut) \, }
 \displaystyle \frac{\Delta \Rightarrow \alpha  \hspace{0.5cm} \Delta, \alpha \Rightarrow \beta}{\Delta \Rightarrow \beta} $$

\noindent {\bf Logic rules}

$$ \mbox{($\wedge \Rightarrow$) \, } \displaystyle \frac{\Delta, \alpha, \beta \Rightarrow \gamma} {\Delta, \alpha \wedge \beta \Rightarrow \gamma} \hspace{2cm} \mbox{($\Rightarrow \wedge$) \, }
 \displaystyle \frac{\Delta \Rightarrow \alpha  \hspace{0.5cm} \Delta \Rightarrow \beta}{\Delta \Rightarrow \alpha \wedge \beta} $$

$$ \mbox{($\vee \Rightarrow$) \, } \displaystyle \frac{\Delta, \alpha \Rightarrow \gamma \hspace{0.5cm} \Delta, \beta \Rightarrow \gamma} {\Delta, \alpha \vee \beta \Rightarrow \gamma}$$

$$ \mbox{($\Rightarrow \vee$) \, }
 \displaystyle \frac{\Delta \Rightarrow \alpha }{\Delta \Rightarrow \alpha \vee \beta} \hspace{2cm} \displaystyle \frac{\Delta \Rightarrow \beta }{\Delta \Rightarrow \alpha \vee \beta}$$

$$ \mbox{($\neg$) \, } \displaystyle \frac{\alpha \Rightarrow \beta } {\neg \beta \Rightarrow \neg \alpha }
 \hspace{2cm}  \mbox{($\bot$) \,} \frac{\Delta \Rightarrow \bot } {\Delta \Rightarrow \alpha}$$

$$ \mbox{($\neg \neg \Rightarrow$) \, } \displaystyle \frac{\Delta, \alpha \Rightarrow \beta } {\Delta, \neg \neg \alpha \Rightarrow \beta} \hspace{2cm} \mbox{($\Rightarrow \neg \neg$)} \, \frac{\Delta \Rightarrow \alpha } {\Delta \Rightarrow \neg \neg \alpha} $$

$$ \mbox{($\square \Rightarrow$) \, } \displaystyle \frac{\Delta, \alpha, \neg \alpha \Rightarrow \beta } {\Delta, \alpha, \neg \square \alpha \Rightarrow \beta} \hspace{2cm} \mbox{($\Rightarrow \square$)} \, \frac{\Delta \Rightarrow \alpha \wedge \neg \alpha } {\Delta \Rightarrow \alpha \wedge \neg \square \alpha} $$

\

\noindent The notion of derivation in  the sequent calculus $\mathfrak{G}$ is the usual. Besides, for every finite set $\Gamma \cup
\{\varphi\} \subseteq Fm$, we write $\Gamma \vdash_{\mathfrak{G}}  \varphi$ iff the sequent $\Gamma \Rightarrow \varphi$ has a derivation in $\mathfrak{G}$. We say that the sequent $\Gamma \Rightarrow \varphi$ is provable iff there exists a derivation for it in $\mathfrak{G}$.\\[2mm]
In~\cite{FR2}, it was proved  that  $\mathfrak{G}$ is sound and complete with respect to the tetravalent modal logic $\cal TML$, constituting therefore a proof-theoretic counterpart of it.

\

\begin{theo}\label{compTML} {\rm(Soundness and Completeness, \cite{FR2})} For every finite set $\Gamma
\cup\{\alpha\} \subseteq Fm$,
$$\Gamma \models_{\cal TML} \alpha \ \ \textrm{ if and only if} \ \ \Gamma \vdash_{\mathfrak{G}} \alpha.$$
\end{theo}
\noindent Moreover,
\begin{prop}{\rm(\cite{FR2})}
An arbitrary equation $\psi \approx \varphi$ holds in every TMA iff $\psi \dashv\vdash_{\mathfrak{G}}\varphi$ (that is, $\psi \vdash_{\mathfrak{G}}\varphi$ and $\varphi \vdash_{\mathfrak{G}}\psi$).
\end{prop}

\noindent As a consequence of it we have that:

\begin{cor}{\rm(\cite{FR2})}\label{CorTheoiffValid}
\begin{itemize}
\item[{\rm(i)}] The equation $\psi \approx \1$ holds in every TMA iff \, $\vdash_{\mathfrak{G}}\psi$.
\item[{\rm(ii)}] For any $\psi, \varphi \in Fm$,\,  $\psi \vdash_{\mathfrak{G}} \varphi$ \, iff \,
 $h(\psi) \leq h(\varphi)$ \, for every \, $h \in \Hom(\mathfrak{Fm}, \mathfrak{A})$, \, for every
 $\mathfrak{A} \in {\bf TMA}$.
\end{itemize}
\end{cor}
\noindent Corollary~\ref{CorTheoiffValid} is a powerful tool to determine whether a given sequent of $\mathfrak{G}$ is provable or not. For instance,
\begin{prop}{\rm(\cite{MF})}\label{prop1} In $\mathfrak{G}$ we have that the sequent $\neg \square \alpha \Rightarrow \alpha$ is provable iff \, the sequent $\Rightarrow \alpha$ is provable.
\end{prop}

\

%\vspace{-6mm}
\noindent Recall that a rule of inference is {\em admissible} in a
formal system if the set of theorems of the system is closed
under the rule; and a rule is said to be {\em derivable} in the
same formal system if its conclusion can be derived from its
premises using the other rules of the system.

\noindent A well--known rule for those readers familiar with modal logic is the {\em Rule of Necessitation},
which states that if $\varphi$ is a theorem, so is $\square \varphi$. Formally,

$$ \mbox{(Nec) \,} \frac{ \Rightarrow \varphi } { \Rightarrow \square \varphi} $$

\noindent Then, we have that:

\begin{lem}{\rm(\cite{MF})}\label{LemNec} The Rule of Necessitation is admissible in $\mathfrak{G}$.
\end{lem}

\begin{prop}{\rm(\cite{MF})}  Every proof of \, $\Rightarrow \square (\alpha \vee \neg \square \alpha)$  in  $\mathfrak{G}$ uses the cut rule.
\end{prop}
\noindent Moreover, we have that, for every $\varphi \in Fm$ such that $\Rightarrow \varphi$  is provable in $\mathfrak{G}$ then $\Rightarrow \square \varphi$ is provable in $\mathfrak{G}$; and every proof of $\Rightarrow \square \varphi$ in $\mathfrak{G}$ makes use of the cut rule (see \cite{MF}). Consequently,

\begin{theo}{\rm(\cite{MF})} $\mathfrak{G}$ does not admit cut--elimination.
\end{theo}

In \cite{CF}, taking profit of the contrapositive implication $\succ$ introduced by A. V. Figallo and P. Landini  in \cite{AVF3}, we introduced sound and complete Hilbert-style calculus for $\tml$. Later,  using a general method proposed by Avron, Ben-Naim and Konikowska (\cite{Avron02}), it was provided a sequent calculus for $\cal TML$ with the cut--elimination property; and, inspired by the latter, it was  presented a {\em natural deduction} system, sound and complete with respect to $\tml$ (\cite{MF}).

\section{Normal proofs for $\tml$}\label{s2}

In this section, we shall present a natural deduction system for ${\cal TML}$. We take our inspiration from the construction  made before. In particular, it threw some light on how the connective $\square$ behaves. The proof system ${\bf ND}_{\cal TML}$ will be defined following the notational conventions given in \cite{TS}.

\begin{defi}\label{defND} Deductions in ${\bf ND}_{\cal TML}$ are inductively defined as follows:\\
Basis: The proof tree with a single occurrence of an assumption $\phi$ with a marker is a deduction with conclusion $\phi$ from open assumption $\phi$ .\\[2mm]
Inductive step: Let ${\cal D}$, ${\cal D}_1$ ,${\cal D}_2$,${\cal D}_3$ be derivations. Then, they can be extended by one of the following rules below. The classes {\rm[$\neg\phi$]$^u$}, {\rm[$\neg\psi$]$^v$}, {\rm[$\phi$]$^u$} , {\rm[$\psi$]$^v$} below contain open assumptions of the deductions of the premises of the final inference, but are closed in the whole deduction. 

\begin{prooftree}
%\alwaysNoLine
%\LeftLabel{\small(Hyp.)}
%\noLine
\AxiomC{}
\RightLabel{\rm MA}
\UnaryInfC{$\phi \vee \neg \square\phi$}
\end{prooftree}

\begin{prooftree}
%\alwaysNoLine
%\LeftLabel{\small(Hyp.)}
%\noLine
\AxiomC{${\cal D}_1$}
\noLine
\UnaryInfC{$\phi$}

\AxiomC{${\cal D}_2$}
\noLine
\UnaryInfC{$\psi$}
\RightLabel{\rm $\wedge$I}
\BinaryInfC{$\phi \wedge \psi$}

\AxiomC{}

\AxiomC{${\cal D}$}
\noLine
\UnaryInfC{$\phi\wedge\psi$}
\RightLabel{\rm $\wedge$E$_1$}
\UnaryInfC{$\phi$}

\AxiomC{}

\AxiomC{${\cal D}$}
\noLine
\UnaryInfC{$\phi\wedge\psi$}
\RightLabel{\rm $\wedge$E$_2$}
\UnaryInfC{$\psi$}

\noLine
\QuinaryInfC{}
\end{prooftree}

\begin{prooftree}
%\alwaysNoLine
%\LeftLabel{\small(Hyp.)}
%\noLine

\AxiomC{${\cal D}$}
\noLine
\UnaryInfC{$\neg\phi$}
\RightLabel{\rm $\neg\wedge$I$_1$}
\UnaryInfC{$\neg(\phi\wedge\psi)$}

\AxiomC{}

\AxiomC{${\cal D}$}
\noLine
\UnaryInfC{$\neg\psi$}
\RightLabel{\rm $\neg\wedge$I$_2$}
\UnaryInfC{$\neg(\phi\wedge\psi)$}

\AxiomC{}

\AxiomC{}
\noLine
\UnaryInfC{${\cal D}_1$}
\noLine
\UnaryInfC{$\neg(\phi\wedge\psi)$}

\AxiomC{\rm[$\neg\phi$]$^u$}
\noLine
\UnaryInfC{${\cal D}_2$}
\noLine
\UnaryInfC{$\chi$}

\AxiomC{\rm [$\neg\psi$]$^v$}
\noLine
\UnaryInfC{${\cal D}_3$}
\noLine
\UnaryInfC{$\chi$}
\RightLabel{\rm $\neg\wedge$E,$u$,$v$}
\TrinaryInfC{$\chi$}
\noLine
\QuinaryInfC{}
\end{prooftree}

\begin{prooftree}
%\alwaysNoLine
%\LeftLabel{\small(Hyp.)}
%\noLine

\AxiomC{${\cal D}$}
\noLine
\UnaryInfC{$\phi$}
\RightLabel{\rm $\vee$I$_1$}
\UnaryInfC{$\phi\vee\psi$}

\AxiomC{}

\AxiomC{${\cal D}$}
\noLine
\UnaryInfC{$\psi$}
\RightLabel{\rm $\vee$I$_2$}
\UnaryInfC{$\phi\vee\psi$}

\AxiomC{}

\AxiomC{}
\noLine
\UnaryInfC{${\cal D}_1$}
\noLine
\UnaryInfC{$\phi\vee\psi$}

\AxiomC{\rm[$\phi$]$^u$}
\noLine
\UnaryInfC{${\cal D}_2$}
\noLine
\UnaryInfC{$\chi$}

\AxiomC{\rm [$\psi$]$^v$}
\noLine
\UnaryInfC{${\cal D}_3$}
\noLine
\UnaryInfC{$\chi$}
\RightLabel{\rm $\vee$E,$u$,$v$}
\TrinaryInfC{$\chi$}
\noLine
\QuinaryInfC{}
\end{prooftree}

\begin{prooftree}
%\alwaysNoLine
%\LeftLabel{\small(Hyp.)}
%\noLine
\AxiomC{${\cal D}_1$}
\noLine
\UnaryInfC{$\neg\phi$}

\AxiomC{${\cal D}_2$}
\noLine
\UnaryInfC{$\neg\psi$}
\RightLabel{\rm $\neg\vee$I}
\BinaryInfC{$\neg(\phi \vee \psi)$}

\AxiomC{}

\AxiomC{${\cal D}$}
\noLine
\UnaryInfC{$\neg(\phi\vee\psi)$}
\RightLabel{\rm $\neg\vee$E$_1$}
\UnaryInfC{$\neg\phi$}

\AxiomC{}

\AxiomC{${\cal D}$}
\noLine
\UnaryInfC{$\neg(\phi\vee\psi)$}
\RightLabel{\rm $\neg\vee$E$_2$}
\UnaryInfC{$\neg\psi$}

\noLine
\QuinaryInfC{}
\end{prooftree}

\begin{prooftree}
%\alwaysNoLine
%\LeftLabel{\small(Hyp.)}
%\noLine
\AxiomC{${\cal D}$}
\noLine
\UnaryInfC{$\phi$}
\RightLabel{\rm $\neg\neg$I}
\UnaryInfC{$\neg\neg\phi$}

\AxiomC{}

\AxiomC{}

\AxiomC{}

\AxiomC{${\cal D}$}
\noLine
\UnaryInfC{$\neg\neg\phi$}
\RightLabel{\rm $\neg\neg$E}
\UnaryInfC{$\phi$}

\noLine
\QuinaryInfC{}
\end{prooftree}

\begin{prooftree}
%\alwaysNoLine
%\LeftLabel{\small(Hyp.)}
%\noLine
\AxiomC{${\cal D}_1$}
\noLine
\UnaryInfC{$\phi$}

\AxiomC{\rm [$\neg\phi$]$^u$}
\noLine
\UnaryInfC{${\cal D}_2$}
\noLine
\UnaryInfC{$\bot$}
\RightLabel{\rm $\square$I,$u$}
\BinaryInfC{$\square\phi$}

\AxiomC{}

\AxiomC{}

\AxiomC{}

\AxiomC{${\cal D}$}
\noLine
\UnaryInfC{$\square\phi$}
\RightLabel{\rm $\square$E}
\UnaryInfC{$\phi$}

\noLine
\QuinaryInfC{}
\end{prooftree}

\begin{prooftree}
%\alwaysNoLine
%\LeftLabel{\small(Hyp.)}
%\noLine
\AxiomC{${\cal D}$}
\noLine
\UnaryInfC{$\neg\phi$}
\RightLabel{\rm $\neg\square$I}
\UnaryInfC{$\neg\square\phi$}

\AxiomC{}

\AxiomC{}

\AxiomC{}

\AxiomC{${\cal D}_1$}
\noLine
\UnaryInfC{$\neg\square\phi$}

\AxiomC{${\cal D}_2$}
\noLine

\UnaryInfC{$\phi$}

\RightLabel{\rm $\neg\square$E}
\BinaryInfC{$\neg\phi$}

\noLine
\QuinaryInfC{}
\end{prooftree}

\begin{prooftree}
%\alwaysNoLine
%\LeftLabel{\small(Hyp.)}
%\noLine

\AxiomC{${\cal D}$}
\noLine
\UnaryInfC{$\neg\phi\wedge\square\phi $}
\RightLabel{$\bot$I}
\UnaryInfC{$\bot$}

\AxiomC{}

\AxiomC{}

\AxiomC{}

\AxiomC{${\cal D}$}
\noLine
\UnaryInfC{$\bot $}
\RightLabel{$\bot$E}
\UnaryInfC{$\alpha$}
\noLine
\QuinaryInfC{}
\end{prooftree}

\end{defi}

\

\begin{rem}\begin{itemize}
\item[]
\item[\rm (i)] Actually, in \cite{MF}, the introduction rule for $\square$ is
\begin{prooftree}
\AxiomC{${\cal D}_1$}
\noLine
\UnaryInfC{$\psi\vee\phi$}

\AxiomC{\rm [$\neg\phi$]$^u$}
\noLine
\UnaryInfC{${\cal D}_2$}
\noLine
\UnaryInfC{$\psi$}
\RightLabel{\rm $\square$I$^*$,$u$}
\BinaryInfC{$\psi\vee\square\phi$}
\end{prooftree}
If we take $\psi$ as $\bot$ in  $\square$I$^*$ we get $\square$I as in Definition \ref{defND}. 
\item[\rm (ii)] The intuition behind rule $\square$E is the following: ``if we have a deduction for $\alpha$ and $\neg \alpha$ is not provable, then we have a deduction for $\square\alpha$''. 
\end{itemize}
\end{rem}

\

As usual, by application of the rule $\neg\wedge$E a new proof-tree is formed from ${\cal D}$, ${\cal D}_1$, and ${\cal D}_2$ by adding at the bottom the conclusion $\chi$ while closing the sets {\rm [$\neg\phi$]$^u$} and {\rm [$\neg\psi$]$^u$} of open  assumptions  marked by $u$ and $v$, respectively. Idem for the rules $\wedge$E and $\square$I. Note that we have introduced the symbol $\bot$, it behaves here as an arbitrary unprovable propositional constant whose negation is provable. \\[2mm]
Let $\Gamma \cup \{\alpha\} \subseteq Fm$, we say that the conclusion $\alpha$ is derivable from a set $\Gamma$ of premises, noted $\Gamma\vdash \alpha$, if and only if there is a deduction in ${\bf ND}_{\cal TML}$ of $\alpha$ from $\Gamma$.

\

\begin{lem} Let $\alpha, \beta\in Fm$. Then, \\[2mm]
\begin{tabular}{clcl}
{\bf (i)} &   $(\alpha \vee  \beta)\wedge\gamma \dashv\vdash (\alpha\wedge\gamma)\vee(\beta\wedge\gamma)$, \hspace{.5cm} & {\bf (ii)} & $(\alpha \wedge  \beta)\vee\gamma \dashv\vdash (\alpha\vee\gamma)\wedge(\beta\vee\gamma)$,\\
{\bf (iii)} &   $\neg(\alpha \vee  \beta) \dashv\vdash \neg\alpha \wedge \neg\beta$, \hspace{1cm} & {\bf (iv)} & $\neg(\alpha \wedge  \beta) \dashv\vdash \neg\alpha \vee \neg\beta$,\\
{\bf (v)} & $\neg\neg \alpha\dashv\vdash \alpha$.  &  & \\
\end{tabular}
\end{lem}
\begin{dem} It is consequence of the I-rules and E-rules for $\wedge$, $\vee$ and $\neg$. 
\end{dem}

\

\noindent Besides, 

\begin{lem}\label{lem} Let $\alpha, \beta\in Fm$. Then, \\[2mm]
\begin{tabular}{clcl}
{\bf (i)} &   $\vdash \square(\alpha \vee  \neg \square\alpha) $, \hspace{1cm} & {\bf (vii)} & $\square \alpha \vee  \neg \alpha \dashv\vdash \alpha \vee  \neg \alpha $,\\

{\bf (ii)} & $\neg\square \alpha\wedge \alpha \dashv\vdash  \neg\alpha\wedge \alpha$  , \hspace{1cm} & {\bf (viii)} & $\square \square \alpha\dashv\vdash  \square\alpha$,\\
{\bf (iii)} & $\vdash \square \alpha \vee \neg\square \alpha$,  & {\bf (ix)} & $\square(a\wedge \beta)\dashv\vdash\square \alpha\wedge \square \beta$, \\
{\bf (iv)} & $\square \alpha \wedge \neg \square \alpha \dashv\vdash \bot$, & {\bf (x)} & $\square (a \vee \square \beta) \dashv\vdash \square \alpha \vee \square \beta$,\\
{\bf (v)} & $\square(\square \alpha\wedge \square \beta)\dashv\vdash\square \alpha\wedge \square \beta$, & {\bf (xi)} & $\square \neg \square \alpha \dashv\vdash \neg\square \alpha$ \\
{\bf (vi)} & $\square(\square \alpha\vee \square \beta)\dashv\vdash\square \alpha\vee \square \beta$, & {\bf (xii)} & $ \alpha\wedge \square \neg \alpha \dashv\vdash \bot$,\\
\end{tabular}
\end{lem}
\begin{dem} We shall only prove {\bf (i)}, {\bf (ii)}, {\bf (ix)} and {\bf (xi)}.\\
{\bf (i)}\\
\begin{prooftree}
\AxiomC{}
\RightLabel{\rm MA}
\UnaryInfC{$\alpha \vee \neg\square\alpha$}

\AxiomC{$\neg(\alpha \vee \neg\square\alpha)^u$}
\doubleLine
\UnaryInfC{$\neg\alpha \wedge \neg\neg\square\alpha)$}
\RightLabel{\small $\wedge$E$_1$}
\UnaryInfC{$\neg\alpha$}

\AxiomC{$\neg(\alpha \vee \neg\square\alpha)^u$}
\doubleLine
\UnaryInfC{$\neg\alpha \wedge \neg\neg\square\alpha)$}
\RightLabel{\small $\wedge$E$_2$}
\UnaryInfC{$\neg\neg\square\alpha$}
\RightLabel{\small $\neg\neg$E}
\UnaryInfC{$\square\alpha$}

\RightLabel{\small $\wedge$I}
\BinaryInfC{$\neg\alpha\wedge\square\alpha$}
\RightLabel{\small($\bot$)}
\UnaryInfC{$\bot$}

\RightLabel{\small $\square$I,$u$}
\BinaryInfC{$\square(\alpha \vee \neg\square\alpha)$}
\end{prooftree}

\noindent {\bf (ii)}:

\begin{prooftree}
\AxiomC{$\neg\alpha \wedge \alpha$}
\RightLabel{\small $\wedge$E$_2$}
\UnaryInfC{$\alpha$}

\AxiomC{$\neg\alpha \wedge \alpha$}
\RightLabel{\small $\wedge$E$_1$}
\UnaryInfC{$\neg\alpha$}
\RightLabel{\small $\neg\square$I}
\BinaryInfC{$\neg\square\alpha$}

\AxiomC{$\neg\alpha \wedge \alpha$}
\RightLabel{\small $\wedge$E$_2$}
\UnaryInfC{$\alpha$}

\RightLabel{\small $\wedge$I}
\BinaryInfC{$\neg\square\alpha\wedge\alpha$}
\end{prooftree}

\

\begin{prooftree}
\AxiomC{$\neg\square\alpha \wedge \alpha$}
\RightLabel{\small $\wedge$E$_1$}
\UnaryInfC{$\neg\square\alpha$}

\AxiomC{$\neg\square\alpha \wedge \alpha$}
\RightLabel{\small $\wedge$E$_2$}
\UnaryInfC{$\alpha$}
\RightLabel{\small $\neg\square$I}
\BinaryInfC{$\neg\alpha$}

\AxiomC{$\neg\square\alpha \wedge \alpha$}
\RightLabel{\small $\wedge$E$_2$}
\UnaryInfC{$\alpha$}

\RightLabel{\small $\wedge$I}
\BinaryInfC{$\neg\alpha\wedge\alpha$}
\end{prooftree}

\noindent {\bf (ix)}:

\begin{prooftree}\small
\AxiomC{$\square(\alpha \wedge \beta)$}
%\RightLabel{\small $\square$E}
\UnaryInfC{$\alpha \wedge \beta$}
%\RightLabel{\small $\wedge$E$_1$}
\UnaryInfC{$\alpha$}

\AxiomC{$\square(\alpha \wedge \beta)$}
%\RightLabel{\small $\square$E}
\UnaryInfC{$\alpha \wedge \beta$}
%\RightLabel{\small $\wedge$E$_1$}
\UnaryInfC{$\alpha$}

%\AxiomC{$\square(\alpha \wedge \beta)$}
%\RightLabel{\small $\square$E}
%\UnaryInfC{$\alpha \wedge \beta$}
%\RightLabel{\small $\wedge$E$_1$}
%\UnaryInfC{$\alpha$}

\AxiomC{$\neg \alpha^u$}
%\RightLabel{\small $\neg\square$I}
\UnaryInfC{$\neg\square\alpha$}
\BinaryInfC{$\alpha\wedge\neg\square\alpha$}
%\RightLabel{\small $\bot$}
\UnaryInfC{$\bot$}
\RightLabel{\small $\square$I,$u$}
\BinaryInfC{$\square\alpha$}

\AxiomC{$\square(\alpha \wedge \beta)$}
%\RightLabel{\small $\square$E}
\UnaryInfC{$\alpha \wedge \beta$}
%\RightLabel{\small $\wedge$E$_2$}
\UnaryInfC{$\beta$}

\AxiomC{$\square(\alpha \wedge \beta)$}
%\RightLabel{\small $\square$E}
\UnaryInfC{$\alpha \wedge \beta$}
%\RightLabel{\small $\wedge$E$_2$}
\UnaryInfC{$\beta$}

%\AxiomC{$\square(\alpha \wedge \beta)$}
%\RightLabel{\small $\square$E}
%\UnaryInfC{$\alpha \wedge \beta$}
%\RightLabel{\small $\wedge$E$_2$}
%\UnaryInfC{$\beta$}

\AxiomC{$\neg \beta^v$}
%\RightLabel{\small $\neg\square$I}
\UnaryInfC{$\neg\square\beta$}
\BinaryInfC{$\beta\wedge\neg\square\beta$}
%\RightLabel{\small $\bot$}
\UnaryInfC{$\bot$}
\RightLabel{\small $\square$I,$v$}
\BinaryInfC{$\square\beta$}
\BinaryInfC{$\square\alpha\wedge\square\beta$}

\end{prooftree}

\

\begin{prooftree}\small
\AxiomC{$\square\alpha\wedge\square\beta$}
%\RightLabel{\small $\wedge$E$_1$}
\UnaryInfC{$\square\alpha$}
\UnaryInfC{$\alpha$}

\AxiomC{$\square\alpha\wedge\square\beta$}
%\RightLabel{\small $\wedge$E$_1$}
\UnaryInfC{$\square\beta$}
\UnaryInfC{$\beta$}
\BinaryInfC{$\alpha\wedge\beta$}

\AxiomC{$\neg(\alpha \wedge \beta)^w$}

\AxiomC{$\square\alpha\wedge\square\beta$}
%\RightLabel{\small $\wedge$E$_1$}
\UnaryInfC{$\square\alpha$}
\UnaryInfC{$\alpha$}

%\AxiomC{$\square\alpha\wedge\square\beta$}
%\RightLabel{\small $\wedge$E$_1$}
%\UnaryInfC{$\square\alpha$}
%\UnaryInfC{$\alpha$}

\AxiomC{$\neg\alpha^u$}
\UnaryInfC{$\neg\square\alpha$}
\BinaryInfC{$\alpha\wedge\neg\square\alpha$}
\UnaryInfC{$\bot$}

\AxiomC{$\square\alpha\wedge\square\beta$}
%\RightLabel{\small $\wedge$E$_1$}
\UnaryInfC{$\square\beta$}
\UnaryInfC{$\beta$}

%\AxiomC{$\square\alpha\wedge\square\beta$}
%\RightLabel{\small $\wedge$E$_1$}
%\UnaryInfC{$\square\beta$}
%\UnaryInfC{$\beta$}

\AxiomC{$\neg\beta^u$}
\UnaryInfC{$\neg\square\beta$}
\BinaryInfC{$\beta\wedge\neg\square\beta$}
\UnaryInfC{$\bot$}

\RightLabel{\small $\neg\wedge$E,$u$,$v$}
\TrinaryInfC{$\bot$}
\RightLabel{\small $\square$I,$w$}
\BinaryInfC{$\square(\alpha\wedge\beta)$}
\end{prooftree}

\noindent {\bf (xi)}: By $\square$E we have that $\square\neg\square\alpha\vdash \neg\square\alpha$. For the converse, consider the following deduction

\begin{prooftree}
\AxiomC{$\neg\square\alpha$}

\AxiomC{$\neg\neg\square\alpha^u$}

\UnaryInfC{$\square\alpha$}

\AxiomC{$\neg\square\alpha$}
\BinaryInfC{$\bot$}
\RightLabel{\small $u$}
\BinaryInfC{$\square\neg\square\alpha$}
\end{prooftree}
\end{dem}

\noindent Note that all syntactic proofs displayed in the Lemma \ref{lem} are normal.

\

\begin{theo}{\rm (Soundness and Completeness, \cite{MF}}) Let $\Gamma, \Delta\subseteq Fm$, $\Gamma$ finite. The following conditions are equivalent:
\begin{itemize} 
\item[\rm (i)] the sequent \ $\Gamma\Rightarrow\Delta$ \ is derivable in ${\bf SC}_{\cal TML}$,
\item[\rm (ii)] there is a deduction of the disjunction of the sentences in $\Delta$ from $\Gamma$ in ${\bf ND}_{\cal TML}$.
\end{itemize}
\end{theo}

\

In what follows, the {\em del-rules} (for ``disjunction-elimination-like") of  ${\bf ND}_{\cal TML}$ are $\vee$E and $\neg\wedge$E. As usual,  a segment (of length $n$) in a deduction ${\cal D}$ of  ${\bf ND}_{\cal TML}$ is a sequence $\alpha_1,\dots , \alpha_n$ of consecutive occurrences of a formula $\alpha$ in {\cal D} such that for $1 \leq n, i < n$,
\begin{itemize}
\item[(a)] $\alpha_i$ is a minor premise of a del-rule application in {\cal D}, with conclusion $\alpha_{i+1}$,
\item[(b)]  $\alpha_n$ is not a minor premise of a del-rule application, and 
\item[(c)] $\alpha_1$ is not the conclusion of a del-rule application. 
\end{itemize}
Note that in this work, the complexity (or degree) of a formula $\alpha$  is defined as follows:

\begin{defi}\label{defi01} Let $\alpha$ be a formula. The complexity (degree) of $\alpha$, $c(\alpha)$, is the natural number obtained by
\begin{itemize}
\item[(i)] if $p$ is a propositional variable then $c(p)=0$,
\item[(ii)] $c(\beta\sharp \gamma)=c(\beta)+c(\gamma)+1$ for $\sharp\in \{\vee, \wedge\}$,
\item[(iii)] $c(\neg\alpha)=c(\alpha) + 1$,
\item[(iv)] $c(\square \alpha)=c(\alpha) + 2$.
\end{itemize}
\end{defi}

A formula occurrence which is neither a minor premise nor the conclusion of an application of a del-rule always belongs to a segment of length 1.\\[2mm]
A segment is {\em maximal}, or a cut (segment) if $\alpha_n$ is the major premise of an E-rule, and either $n > 1$, or $n = 1$ and $\alpha_n$ is the conclusion of an I-rule. The {\em cutrank} $cr(\sigma)$ of a maximal segment $ \sigma$ with formula $\alpha$ is the complexity of $\alpha$. %Note that in this work, the complexity (or degree) of a formula $\alpha$ ($c(\alpha)$) is defined as follows: if $p$ is a propositional variable then $c(p)=0$; $c(\beta\sharp \gamma)=c(\beta)+c(\gamma)+1$ for $\sharp\in \{\vee, \wedge\}$; $c(\neg\alpha)=c(\alpha) + 1$; and $c(\square \alpha)=c(\alpha) + 2$.\\
Besides, the cutrank $cr({\cal D})$ of a deduction ${\cal D}$ is the maximum of the cutranks of cuts in ${\cal D}$. If there is no cut, the cutrank of ${\cal D}$ is zero. A {\em critical cut} of ${\cal D}$ is a cut of maximal cutrank among all cuts in ${\cal D}$.  A deduction without critical cuts is said to be {\em normal}.

\begin{lem} Let ${\cal D}$ be a deduction in ${\bf ND}_{\cal TML}$. Then, ${\cal D}$ reduces to a deduction ${\cal D}'$ in which the consequence of every application of the $\bot$E rule is a propositional variable $p$ (atomic) or its negation $\neg p$.  
\end{lem} \begin{dem} Consider the following deduction with an application of the $\bot$E
\begin{center}
\AxiomC{$\cal D$}
\noLine
\UnaryInfC{$\bot$}
\RightLabel{\tiny $\bot$E}
\UnaryInfC{$\alpha$}
\DisplayProof
\end{center}
It is not difficult to check that if $\alpha$ has the shape of $\gamma_1\wedge\gamma_2$, $\gamma_1\vee\gamma_2$, $\neg(\gamma_1\wedge\gamma_2)$, $\neg(\gamma_1\vee\gamma_2)$, $\neg\neg\gamma_1$, $\square\gamma_1$ and $\neg\square\gamma_1$ we can remove this application of $\bot$E using  application(s) of $\bot$E with consequence formula(s) that has complexity strictly less that the complexity of $\alpha$. Thus, by successively repeating this transformation we can finally obtain a deduction with the required characteristics.\end{dem}

\
    
In what follows, we shall only consider deductions in which  the consequence of every application of the $\bot$E rule is a propositional variable $p$ (atomic) or its negation $\neg p$. Observe that, in this kind of deductions, there can not be occurrences of formulas that are consequence of $\bot$E and the major premise of an E-rule.	

\

We first show how to remove cuts of length $1$. Besides, $\wedge$-conversions and $\vee$-conversions are as in the system of natural deduction for intuitionistic (or classical) logic.

\

\noindent $\neg\wedge$-{\em conversion}
\begin{center}
\AxiomC{${\cal D}$}
\noLine
\UnaryInfC{$\neg\alpha_i$}
\UnaryInfC{$\neg(\alpha_1\wedge\alpha_2)$}

\AxiomC{\rm[$\neg\alpha_1$]$^u$}
\noLine
\UnaryInfC{${\cal D}_1$}
\noLine
\UnaryInfC{$\chi$}

\AxiomC{\rm [$\neg\alpha_2$]$^v$}
\noLine
\UnaryInfC{${\cal D}_2$}
\noLine
\UnaryInfC{$\chi$}
\RightLabel{\rm $u$,$v$}
\TrinaryInfC{$\chi$}
\DisplayProof \, \,
converts to \, \,
\AxiomC{${\cal D}$}
\noLine
\UnaryInfC{[$\neg\alpha_i$]}
\noLine
\UnaryInfC{${\cal D}_i$}
\noLine
\UnaryInfC{$\chi$}
\DisplayProof for $i=1,2.$\end{center}
\noindent $\neg\vee$-{\em conversion}
\begin{center}
\AxiomC{${\cal D}_1$}
\noLine
\UnaryInfC{$\neg\alpha_1$}
\AxiomC{${\cal D}_2$}
\noLine
\UnaryInfC{$\neg\alpha_2$}
\BinaryInfC{$\neg(\alpha_1\wedge\alpha_2)$}
\UnaryInfC{$\neg\alpha_i$}
\DisplayProof \, \,
converts to \, \,
\AxiomC{${\cal D}_i$}
\noLine
\UnaryInfC{$\neg\alpha_i$}
\DisplayProof for $i=1,2.$\\[2mm]
\end{center}
\noindent $\neg\neg$-{\em conversion}
\begin{center}
\AxiomC{${\cal D}$}
\noLine
\UnaryInfC{$\alpha$}
\UnaryInfC{$\neg\neg\alpha$}
\UnaryInfC{$\alpha$}
\DisplayProof \, \,
converts to \, \,
\AxiomC{${\cal D}$}
\noLine
\UnaryInfC{$\alpha$}
\DisplayProof 
\end{center}
\noindent $\square$-{\em conversion}
\begin{center}
\AxiomC{${\cal D}_1$}
\noLine
\UnaryInfC{$\alpha$}
\AxiomC{[$\neg\alpha$]$^u$}
\noLine
\UnaryInfC{${\cal D}_2$}
\noLine
\UnaryInfC{$\bot$}
\BinaryInfC{$\square\alpha$}
\UnaryInfC{$\alpha$}
\DisplayProof \, \,
converts to \, \,
\AxiomC{${\cal D}_1$}
\noLine
\UnaryInfC{$\alpha$}
\DisplayProof 
\end{center}
\noindent $\neg\square$-{\em conversion}
\begin{center}
\AxiomC{${\cal D}_1$}
\noLine
\UnaryInfC{$\neg\alpha$}
\UnaryInfC{$\neg\square\alpha$}

\AxiomC{${\cal D}_2$}
\noLine
\UnaryInfC{$\alpha$}
\BinaryInfC{$\neg\alpha$}
\DisplayProof \, \,
converts to \, \,
\AxiomC{${\cal D}_1$}
\noLine
\UnaryInfC{$\neg\alpha$}
\DisplayProof 
\end{center}

\

\noindent In order to remove cuts of length $>1$, we permute E-rules upwards over minor premises of del-rules:\\
\begin{center}
\AxiomC{${\cal D}$}
\noLine
\UnaryInfC{$\alpha_1\vee\alpha_2$}

\AxiomC{\rm[$\alpha_1$]$^u$}
\noLine
\UnaryInfC{${\cal D}_1$}
\noLine
\UnaryInfC{$\chi$}

\AxiomC{\rm [$\alpha_2$]$^v$}
\noLine
\UnaryInfC{${\cal D}_2$}
\noLine
\UnaryInfC{$\chi$}
\RightLabel{\rm $\vee$E}
\TrinaryInfC{$\chi$}
\AxiomC{${\cal D}'$}
\RightLabel{\rm E-rule}
\BinaryInfC{$\chi$'}
\DisplayProof \, \,
converts to \, \,
\AxiomC{${\cal D}$}
\noLine
\UnaryInfC{$\alpha_1\vee\alpha_2$}

\AxiomC{\rm[$\alpha_1$]$^u$}
\noLine
\UnaryInfC{${\cal D}_1$}
\noLine
\UnaryInfC{$\chi$}
\AxiomC{${\cal D}'$}
\RightLabel{\rm E-rule}
\BinaryInfC{$\chi$'}

\AxiomC{\rm [$\alpha_2$]$^v$}
\noLine
\UnaryInfC{${\cal D}_2$}
\noLine
\UnaryInfC{$\chi$}
\AxiomC{${\cal D}'$}
%\RightLabel{\rm E-rule}
\BinaryInfC{$\chi$'}
%\RightLabel{\rm $\vee$E,$u$,$v$}
\TrinaryInfC{$\chi$'}
\DisplayProof,
\end{center}

\begin{center}
\AxiomC{${\cal D}$}
\noLine
\UnaryInfC{$\neg(\alpha_1\wedge\alpha_2)$}

\AxiomC{\rm[$\neg\alpha_1$]$^u$}
\noLine
\UnaryInfC{${\cal D}_1$}
\noLine
\UnaryInfC{$\chi$}

\AxiomC{\rm [$\neg\alpha_2$]$^v$}
\noLine
\UnaryInfC{${\cal D}_2$}
\noLine
\UnaryInfC{$\chi$}
\RightLabel{\rm $\vee$E}
\TrinaryInfC{$\chi$}
\AxiomC{${\cal D}'$}
\RightLabel{\rm E-rule}
\BinaryInfC{$\chi$'}
\DisplayProof \, \,
converts to \, \,
\AxiomC{${\cal D}$}
\noLine
\UnaryInfC{$\neg(\alpha_1\vee\alpha_2)$}

\AxiomC{\rm[$\neg\alpha_1$]$^u$}
\noLine
\UnaryInfC{${\cal D}_1$}
\noLine
\UnaryInfC{$\chi$}
\AxiomC{${\cal D}'$}
\RightLabel{\rm E-rule}
\BinaryInfC{$\chi$'}

\AxiomC{\rm [$\neg\alpha_2$]$^v$}
\noLine
\UnaryInfC{${\cal D}_2$}
\noLine
\UnaryInfC{$\chi$}
\AxiomC{${\cal D}'$}
%\RightLabel{\rm E-rule}
\BinaryInfC{$\chi$'}

%\RightLabel{\rm $\vee$E,$u$,$v$}
\TrinaryInfC{$\chi$'}
\DisplayProof,
\end{center}

\

Applications of $\vee$E ($\neg\wedge$E) with major premise $\alpha_1\vee\alpha_2$ ($\neg(\alpha_1\wedge\alpha_2)$), where at least one of [$\alpha_1$], [$\alpha_2$] ([$\neg\alpha_1$], [$\neg\alpha_2$]) is empty in the deduction of the first or second minor premise, are redundant and can be removed easily. Now, we are in conditions to prove the main result of this section.

\

\begin{theo} (Normalization) Every deduction ${\cal D}$ in ${\bf ND}_{\cal TML}$ reduces to a normal deduction.
\end{theo} 
\begin{dem} We assume that, in every application of an E-rul, the major premise is always to the left of the minor premise(s), if there are any minor premises. We use a main induction on the cutrank $n$ of ${\cal D}$, with a subinduction on $m$, the sum of lengths of all critical cuts (= cut segments) in ${\cal D}$. By a suitable choice of the critical cut to which we apply a conversion we can achieve that either $n$ decreases (and we can appeal to the main induction hypothesis), or that $n$ remains constant but $m$ decreases (and we can appeal to the subinduction hypothesis). Let $\sigma$ be a top critical cut in ${\cal D}$ if no critical cut occurs in a branch of ${\cal D}$ above $\sigma$. Now apply a conversion to the rightmost top critical cut of ${\cal D}$; then the resulting ${\cal D}'$ has a lower cutrank (if the segment treated has length 1, and is the only maximal segment in D), or has the same cutrank, but a lower value for $m$.  
\end{dem}

\

\section{\large \bf The contrapositive implication in $\tml$}
\label{s3}

The original language of the logic of TMAs -- in
particular, the language of logic $\tml$-- does not have an
implication as a primitive connective. It is a natural question
to ask how to define a binary operator in TMAs, in
terms of the others, with the behavior of an implication. Such
operators are useful in order to characterize the lattice of
congruences of a given class of algebras.

Some proposal for an implication operator in TMAs
appeared in the literature. For instance, I. Loureiro proposed in~\cite{L1}
the following implication for TMAs:
$$x \rightarrow y = \neg \square x \vee y;$$
and, by considering the operator
$$x \mapsto y = (x \rightarrow y) \wedge (\neg \square
\neg y \vee \neg x),$$
A. V. Figallo and P. Landini introduced in~\cite{AVF3} an
interesting implication operator for TMAs defined
as follows:
$$x \succ y = (x \mapsto y) \wedge ((\neg x \vee y)\rightarrow(\square  \neg x \vee y)).$$
This operator was called in~\cite{AZ} {\em
contrapositive implication} for TMAs.\\
It is easy to see that the contrapositive implication can be written in a simpler form:
$$x \succ y = (x \rightarrow y) \wedge (\neg y \rightarrow \neg x) \wedge ((\neg x \vee y) \rightarrow (\square  \neg x \vee y)).$$
The main feature of the contrapositive implication is that it internalizes 
the consequence relation (whenever just one premise is considered), as we shall see in Theorem~\ref{MTD}.
Another important aspect of the contrapositive implication is that
all the operations of the TMAs can be defined in terms of $\succ$ and $\0$. In fact:

\begin{prop} [cf.~\cite{AVF3}] \label{DefCon} In every 
TMA it holds:
\begin{multicols}{2}
\begin{itemize}
\item[{\rm(i)}] $\1 \approx (\0 \succ \0)$,
\item[{\rm(ii)}] $\neg x \approx (x \succ \0)$,
\item[{\rm(iii)}] $x \vee y \approx (x \succ y)\succ y$,
\item[{\rm(iv)}] $x \wedge y \approx \neg(\neg x \vee \neg y)$,
\item[{\rm(v)}] $\square x \approx \neg(x \succ \neg x)$.
\end{itemize}
\end{multicols}
Therefore, $\succ$ and $\0$ are enough to generate all the
operations of a given TMA.
\end{prop}

\noindent Additionally, from Proposition~\ref{DefCon} an axiomatization for
 TMAs was given in~\cite{AVF3} in terms of $\succ$
and $\0$ as follows.

\begin{prop} [cf.~\cite{AVF3}] \label{DefTMA}
In every TMA it can be defined a binary operation
$\succ$ and an element $\0$ such that the following holds:
\begin{multicols}{2}
\begin{itemize}
\item[{\rm(C1)}] $(\1 \succ x) \approx x$,
\item[{\rm(C2)}] $(x \succ \1) \approx \1$,
\item[{\rm(C3)}] $(x \succ y) \succ y \approx (y \succ x) \succ x$,
\item[{\rm(C4)}] $x \succ (y \succ z)\approx \1$ implies $y \succ (x \succ z)\approx \1$,
\item[{\rm(C5)}] $((x \succ (x \succ y)) \succ x) \succ x \approx \1$,
\item[{\rm(C6)}] $(\0 \succ x) \approx \1$,
\item[{\rm(C7)}]  $(x \succ \0) \approx \neg x$,
\item[{\rm(C8)}] $((x \wedge y) \succ z)\succ((x \succ z) \vee (y \succ z))= \1$.
\end{itemize}
\end{multicols}
\noindent Conversely, if an algebra with a binary operation
$\succ$  and an element $0$ satisfies (C1)-(C8) where
$\1=_{def} \0\succ \0$, $x \vee y =_{def}
(x \succ y)\succ y$ and $x \wedge y =_{def} \neg(\neg x \vee \neg
y)$, then the resulting structure is a TMA where
$\square x =_{def} \neg(x \succ \neg x)$.
\end{prop}

\begin{defi} \label{defTAMc}
A {\em contrapositive tetravalent modal algebra} is an algebra
$\langle A, \succ, 0 \rangle$ of type $(2,0)$ that satisfies
items (C1)--(C6) and (C8) of Proposition~\ref{DefTMA} (with the
abbreviations defined therein). We shall denote the class of
these algebras by ${\bf TMA}^{c}$.
\end{defi}

\

\noindent Observe that the classes {\bf TMA} and {\bf
TMA}$^{c}$ are termwise equivalent. The main differences reside in the
underlying language defining both classes and the fact that the characterization of the
latter does not allow to see that in fact it is a variety. As it was showed in \cite{CF}, the contrapositive implication $\succ$ is very useful when describing a Hilbert-style system for $\tml$. 

It is worth noting that in $\mathfrak{M}_{4m}$, the canonical TMA,
the contrapositive implication $\succ$ has the following truth-table:

\

\begin{center}

%\begin{table}[h]
\begin{tabular}{|c|c|c|c|c|}         \hline
$\succ$ &$\0$ & $\n$ & $\b$ & $\1$\\ \hline
$\0$       & $\1$ & $\1$ & $\1$ & $\1$ \\ \hline
$\n$       & $\n$ & $\1$ & $\b$ & $\1$         \\ \hline
$\b$       & $\b$ & $\n$ & $\1$ & $\1$         \\ \hline
$\1$       & $\0$ & $\n$ & $\b$ & $\1$         \\ \hline
\end{tabular}
\end{center}

\

\begin{rem} \label{M4=TMAc}
Clearly, the logic of the contrapositive tetravalent modal algebras $\models_{TMA^{c}}$ can be defined by analogy with Definition~\ref{logTMA}.
\end{rem}

\noindent The new connective $\succ$ has some nice
properties, displayed below:

\begin{prop}\label{PropModales} Let $\alpha, \beta \in Fm$.
Then the following holds in $\tml$:\\
\begin{multicols}{2}
\begin{enumerate}[\rm (i)]
\item $\models_{\tml}\bot \succ \alpha$,
\item $\models_{\tml}\alpha \succ \top$,
\item $\models_{\tml} \alpha \succ (\beta \succ \alpha)$,
\item $\models_{\tml} (\alpha \vee \beta) \succ (\beta \vee \alpha)$,
\item $\models_{\tml} \neg \neg \alpha \succ \alpha$,
\item $\models_{\tml} \alpha \succ \neg \neg \alpha$,
\item $\models_{\tml} \square \alpha \succ \square\square \alpha$,
\item $\models_{\tml} \square\alpha \succ \alpha$,
\item $\models_{\tml} \alpha \succ \square \diamondsuit \alpha$,
\item $\models_{\tml} \square \alpha \succ \diamondsuit \alpha$,
\end{enumerate}
\end{multicols}
\hspace*{-0.4cm}{\rm(xi)} $\models_{\tml} \square (\alpha \succ \beta) \succ (\square \alpha \succ \square
\beta)$,\\[2mm]
\, {\rm(xii)} $\models_{\tml} (\diamondsuit \alpha \wedge \diamondsuit \beta)\succ (\diamondsuit ( \alpha \wedge \diamondsuit \beta) \vee \diamondsuit ( \alpha \wedge \beta) \vee \diamondsuit ( \beta\wedge \diamondsuit \alpha)) $.
\end{prop}
\begin{dem} By the simple inspection of the truth-tables.
\end{dem}

\

\begin{rem}\label{Rem2} Theorem (xi) is the {\bf (K)} axiom  {\rm(}see
{\rm\cite{ModLog})}. Theorems (vii), (viii), (ix), (x) and (xii)
correspond to the axioms {\bf (4)}, {\bf (T)}, {\bf (B)}, {\bf
(D)} and {\bf (.3)}, respectively {\rm(}see {\rm\cite{ModLog})}.
Therefore $\tml$ satisfies all the modal axioms of
 the classical modal logic {\bf S5}.  Nevertheless, we cannot affirm that
$\tml$ is a normal modal logic since the implication $\succ$
does not satisfy some properties of the classical implication
{\rm(}see
{\rm \cite{ModLog})}. 

There exist interesting similarities between \L ukasiewicz's
\L$_3$ (seen as a modal logic) and $\tml$. In both logics,
$\square\alpha$ and $\diamondsuit\alpha$ are defined by the same
formulas, namely $\neg(\alpha \succ\neg\alpha)$ and $\neg\alpha
\succ\alpha$, respectively (in the case of \L$_3$, $\neg$ and
$\succ$ denote the respective negation and implication
operators). Moreover, both implications (\L$_3$'s implication and the
contrapositive implication) do not satisfy the contraction law:
$\alpha \succ(\alpha \succ \beta)$ is not equivalent to $(\alpha
\succ \beta)$. From this, both logics satisfy the following modal
principle: $\alpha \succ(\alpha \succ \square\alpha)$, which is
not valid in the classical modal logic {\bf S5}.
\end{rem}

\

\noindent Let $\square^{0}\alpha=_{def} \alpha$ and
$\square^{n+1}\alpha=_{def} \square^{n} \square \alpha$ for any
$n \in \mathbb{N}$.  Analogously, it is defined $\diamondsuit^{n}
\alpha$. Then,

\begin{prop} $\tml$ satisfies the following well-known instance of the Lemmon-Scott schemes {\rm(cf.~\cite{LeSc})} for any $n,l,k,m \in
\mathbb{N}$,

$$\models_{\tml} \diamondsuit^{k} \square^{l} \alpha \succ \square^{m} \diamondsuit^{n} \alpha,$$

\noindent but $\not \models_{\tml} \square \diamondsuit \alpha
\succ \diamondsuit \square \alpha$.
\end{prop}
\begin{dem} From Proposition~\ref{PropModales}.
\end{dem}

\

\noindent Finally, in $\tml$ we have a weak version of the {\em
Deduction Metatheorem} with respect to the contrapositive
implication.

\begin{theo} [\cite{AVF3}] \label{MTD} Let $\alpha, \beta \in Fm$. The following conditions are equivalent:
\begin{itemize}
\item[{\rm(i)}] $\alpha \models_{\tml} \beta$,
\item[{\rm(ii)}] $\models_{\tml} \alpha \succ \beta$.
\end{itemize}
\end{theo}

\

\noindent This last result shows that the contrapositive
implication $\succ$ internalizes the consequence relation of
$\tml$ whenever just one premise is considered. In algebraic terms, $\succ$ internalizes the partial
order $\leq$ of TMAs.\\
It is worth noting that it is not possible to improve Theorem~\ref{MTD} in the following sense:

\begin{prop} \label{not-MTD} In $\tml$ both directions of the deduction metatheorem, with respect to $\succ$,  
fail if more than one premise are allowed. Specifically:
\begin{itemize}
\item[{\rm(i)}] $\alpha, \beta \models_{\tml} \gamma$ \ does not imply that\ $\alpha \models_{\tml} \beta \succ \gamma$,
\item[{\rm(ii)}] $\alpha \models_{\tml} \beta \succ \gamma$ \ does not
imply that \  $\alpha, \beta \models_{\tml} \gamma$.
\end{itemize}
\end{prop}

\begin{dem}
(i) Observe that $\n \wedge \b \leq \0$, but $\n \not\leq \b \succ
\0=\b$. In order to find an example of this, consider
$\alpha=(\bullet p \wedge \bullet q \wedge \bullet (p \succ q)
\wedge p)$, $\beta=q$ and $\gamma=\bot$, where $p,q$ are two
different propositional variables, and $\bullet\delta=_{def}
\diamondsuit (\delta \wedge \neg \delta)$ is the {\em
inconsistency} operator. Then
$h(\alpha \wedge \beta)=\0=h(\gamma)$ for every $h \in
\Hom(\mathfrak{Fm}, \mathfrak{M}_{4m})$. That is, $\alpha, \beta
\models_{\tml} \gamma$. Now, let $h$ such that $h(p)=\n$ and
$h(q)=\b$. Then $h(\alpha)=\n$ and $h(\beta)=\b$ and so
$h(\alpha)\not\leq h(\beta \succ \gamma)= \b \succ \0=\b$. Therefore, $\alpha \not\models_{\tml} \beta \succ \gamma$. \\[2mm]
(ii) Note that $\n \leq \n \succ \0=\n$, but $\n \wedge \n = \n \not\leq
\0$. For instance, consider $\alpha=p$, $\beta=\neg p$ and
$\gamma=\bot$, where $p$ is a propositional variable. Then $\alpha
\models_{\tml} \beta \succ \gamma$, since
$\alpha\models_{\tml}\neg\neg\alpha$. Consider $h \in
\Hom(\mathfrak{Fm}, \mathfrak{M}_{4m})$ such that $h(p)=\n$. Then
$h(\alpha \wedge \beta)=\n \not\leq \0=h(\gamma)$ and so $\alpha,
\beta \not\models_{\tml} \gamma$.
\end{dem}

\

\noindent In particular, the contrapositive implication does not
satisfy local \emph{modus ponens} (in the sense of~\cite{CCGGS}).

\noindent The importance of the contrapositive implication was confirmed in \cite{CF}, where a Hilbert-style
axiomatization for $\tml$ was given in terms of $\succ$ and $\bot$.

\section{A tableau system for $\tml$ in the language $\{\neg, \succ\}$}\label{s4}

In this section, by adapting the general
techniques introduced in~\cite{CCCM}, we define a decidable
tableau system adequate to check validity in the logic $\tml$.
This constitutes a new (proof-theoretic) decision procedure 
for checking validity in the variety of tetravalent modal
algebras, besides the four-valued truth-tables of $\tml$ and the one available in terms of the cut-free sequent system introduced in \cite{MF}.

The procedure for finding a set of tableau-rules for $\tml$ is
based on the general method presented in~\cite{CCCM} for obtaining
bivalued semantics and  tableau rules for a wide class of finite
matrix logics (see~\cite{CalMar} for further development of this
technique). The given matrix logic must satisfy just one
condition: to be expressive enough to ``separate'' among the
different truth-values of the same kind, namely distinguished and
non-distinguished.

Since we are here interested in just one example, the logic
$\tml$, we will simplify the procedure for obtaining the tableau
rules without entering in the (rather involved) details of the
general construction presented in~\cite{CCCM} and~\cite{CalMar}. Moreover, for the sake of the reader, we will present original proofs of soundness and completeness of the
generated tableau system, generalizing the classical ones of~\cite{Smu}, and so this section will be self-contained.

For simplicity, and taking profit of the contrapositive implication, we will use the
language $\mathfrak{Fm}''=\langle Fm'', \succ, \neg\rangle$. The use of $\neg$ instead of $\bot$ as primitive will be convenient in order to simplify the presentation of the rules, besides the fact that the negation $\neg$ will play a fundamental role in the sequel. It is worth noting that the ``expressive power'' of $\mathfrak{Fm}''$ is the same than that of $\mathfrak{Fm}'$, since $\bot$ can be defined in the former as $\neg(\alpha \succ\alpha)$ for any $\alpha$. Additionally, the tableau rules will be extended to the usual language $\mathfrak{Fm}$ of tetravalent modal algebras.

In the following subsections, we will assume that the reader is
acquainted with the definition of tableau, as well as with the
related notions such as closure rules, closed and open branches
etc. The reader unfamiliar with such concepts can consult the
classical book~\cite{Smu}.

\

\subsection{Separating the truth-values of $\tml$}

From now, $\tml$ will be seen as the matrix logic ${\cal
M}_N=\langle\mathfrak{M}_{4m},\{\b, 1\}\rangle$. Consider the function \break $f:M_4 \to
\{\T,\F\}$ given by $f(\1)=f(\b)=\T$ and $f(\0)=f(\n)=\F$. This function
splits the truth-values into two classes: the distinguished and
the non-distinguished ones.

\noindent Consider now the formula $\neg p$ in $Fm''$. This formula (seen as an
operator over $M_4$) ``separates'' the truth-values of $\tml$ as follows: given $x \in M_4$,\\
$$x = \1 \, \mbox{ iff } \, \, f(x)=\T\mbox{ and } f(\neg x)=\F;$$
$$x = \b \, \mbox{ iff } \, \, f(x)=\T \mbox{ and } f(\neg x)=\T;$$
$$x = \n \, \mbox{ iff } \, \, f(x)=\F\mbox{ and } f(\neg x)=\F;$$
$$x = \0 \, \mbox{ iff } \, \, f(x)=\F \mbox{ and } f(\neg x)=\T.$$
From this it follows:

\begin{center}
$(\ddagger) \left \{ \begin{tabular}{rl}
$x \in \{\1,\b\}$  & iff  $f(x)=\T$;\\
$x \in \{\1,\n\}$ & iff $f(\neg x)=\F$;\\
$x \in \{\0,\b\}$ & iff $f(\neg x)=\T$;\\
$x \in \{\0,\n\}$ & iff $f( x)=\F$.\\
\end{tabular}\right.$\\
\end{center}

\

\subsection{Describing the truth-table of $\succ$ in terms of $\T/\F$}

\noindent By inspection of the truth-table of $\succ$ and by
using $(\ddagger)$ it follows, for every $x,y \in M_4$:\\[2mm]

$f(x\succ y) = \T\ \textrm{ iff } \ \left \{ \begin{tabular}{lll}
 & $f(y)=\T$ & OR\\[2mm]
$f(\neg x)=\T$, & $f(y)= \F$, $f(\neg
y)=\T$ & OR\\[2mm]
$f(x)=\F$ & $f(y)=\F$, $f(\neg y)=\F$ & \\[2mm]
\end{tabular}\right.$\\[2mm]

\

$f(x\succ y) = \F \ \textrm{ iff } \ \left \{ \begin{tabular}{lll}
$f(x)=\T$, & $f(y)= \F$, $f(\neg
y)=\F$, & OR\\[2mm]
$f(\neg x)=\F$,  & $f(y)=\F$, $f(\neg y)=\T$. & \\[2mm]
\end{tabular}\right.$\\[2mm]

\

$f(\neg(x\succ y)) = \T\textrm{ iff } \left \{
\begin{tabular}{lll} $f(x)=\T$, &
$f(y)= \F$, $f(\neg y)=\T$ & OR\\[2mm]
$f(\neg x)=\F$, & $f(y)= \T$, $f(\neg
y)=\T$. & \\[2mm]
\end{tabular}\right.$\\[2mm]

\

$f(\neg(x\succ y)) = \F \textrm{ iff } \left \{
\begin{tabular}{lll}  & $f(\neg y)=\F$ & OR\\[2mm]
$f(\neg x)=\T$, & $f(y)= \T$, $f(\neg
y)=\T$ & OR\\[2mm]
$f(x)=\F$, & $f(y)=\F$, $f(\neg y)=\T$ & \\[2mm]
\end{tabular}\right.$\\[2mm]

\

\subsection{Obtaining the tableau rules for $\tml$}

By substituting in the expressions above the truth-values $x$,
$y$ by formulas $\alpha$, $\beta$ of $Fm''$, and by substituting
equations ``$f(x)=\T$'' and ``$f(x)=\F$'' by signed formulas
$T(\alpha)$ and $F(\alpha)$, respectively, we obtain
automatically the following tableau rules for $\tml$:\\

$$\frac{T(\alpha \succ\beta)}{ T(\beta) \,  |  \,  T(\neg\alpha), F(\beta), T(\neg\beta) \, | \,  F(\alpha), F(\beta), F(\neg\beta)}$$

\

$$\begin{array}{c}
 \displaystyle \frac{F(\alpha
\succ\beta)}{T(\alpha), F(\beta), F(\neg\beta)\ | \ F(\neg\alpha),
F(\beta), T(\neg\beta)}
\end{array}$$

\

$$\begin{array}{c}
 \displaystyle \frac{T(\neg(\alpha
\succ\beta))}{T(\alpha), F(\beta), T(\neg\beta)\ | \
F(\neg\alpha), T(\beta), T(\neg\beta)}
\end{array}$$

\

$$\frac{F(\neg(\alpha \succ\beta))}{F(\neg\beta) \,  | \, T(\neg \alpha),
T(\beta), T(\neg\beta) \, | \, F(\alpha), F(\beta), T(\neg\beta)}$$

\

$$\begin{array}{ccc}
\displaystyle \frac{T(\neg\neg\alpha)}{T(\alpha)} &
\mbox{\ \ \ \ } & \displaystyle
\frac{F(\neg\neg\alpha)}{F(\alpha)}
\end{array}$$

\

\noindent {\bf Closure rule:}
$$
\displaystyle \frac{T(\alpha), F(\alpha)}{\star}$$

\

\noindent Let $\mathbb{T}$ be the tableau system defined by the
rules above. Given a signed formula $\eta$, a completed tableau
starting from $\eta$ is called {\em a tableau for $\eta$}. We say
that a formula $\alpha$ in $Fm''$ is {\em provable in
$\mathbb{T}$}, written as $\vdash_{\mathbb{T}} \alpha$, if there
exists a closed tableau for the signed formula $F(\alpha)$.

\

Given a tableau system, it is in general convenient to define
derived rules in order to get shorter proofs. So, we state a fundamental set of derived rules:

\

\begin{prop} \label{der-rule}
The following rules can be derived in  $\mathbb{T}$:

$$\frac{T(\alpha \succ\beta), T(\neg(\alpha \succ\beta))}
{ F(\neg \alpha), T(\beta), T(\neg\beta) \, | \, T(\alpha), T(\neg\alpha), F(\beta),  T(\neg\beta)} \hspace{1cm}  
\frac{F(\alpha \succ\beta), F(\neg(\alpha \succ\beta))}
{T(\alpha), F(\beta), F(\neg\beta) \, | \, F(\alpha), F(\neg\alpha), F(\beta), T(\neg\beta)}
$$

$$
\frac{T(\alpha \succ\beta), F(\neg(\alpha \succ\beta))}{ F(\alpha), T(\neg\alpha) \, | \, F(\alpha), F(\neg\alpha), F(\beta), F(\neg \beta) \, | \, T(\alpha), T(\neg\alpha), T(\beta), T(\neg \beta) \, |\, T(\beta), F(\neg\beta)}
$$

$$\frac{F(\alpha \succ\beta), T(\neg(\alpha \succ\beta))}{T(\alpha), F(\neg\alpha), F(\beta), T(\neg\beta)}
$$

\

\end{prop}
\begin{dem}
Straightforward, by using the rules of $\mathbb{T}$.
\end{dem}

\subsection{Soundness and completeness of  $\mathbb{T}$}\label{s5.4}

Now we will prove the adequacy of $\mathbb{T}$, that is, that
$\vdash_{\mathbb{T}} \alpha$ if and only if $\models_{\tml}
\alpha$, for every formula $\alpha$. We begin by introducing some
definitions and previous results.

\noindent Given a formula $\alpha \in Fm''$, the {\em degree} of
$\alpha$, denoted by $d(\alpha)$, is a natural number defined as
follows: $d(p)=1$ (for $p \in Var$); $d(\alpha \succ
\beta) = d(\alpha) + d(\beta) + 1$;
$d(\neg\alpha)= d(\alpha)+1$.

\noindent It is clear that the degree of the formulas occurring
in the conclusion of a rule of $\mathbb{T}$ is strictly
less than the degree of the premise of the rule. From this, it is
immediate to prove the following useful result:

\begin{prop} \label{completabl}
Given a signed formula $\eta$, it is always possible to build a
(open or closed) completed tableau for $\eta$.
\end{prop}

\noindent Given an homomorphism $h: \mathfrak{Fm}'' \to
\mathfrak{M}_{4m}$ and a signed formula $\eta$, we say that {\em
$h$ satisfies $\eta$} if
\begin{description}
  \item[--] $\eta=T(\alpha)$  and $h(\alpha) \in \{\1, \b\}$;
  \item[--] $\eta=F(\alpha)$  and $h(\alpha) \in \{\0, \n\}$;
\end{description}

\noindent It follows that $h$ satisfies $T(\neg\alpha)$  iff $h(\alpha) \in \{\0,\b\}$, and $h$ satisfies $F(\neg\alpha)$ iff $h(\alpha) \in \{\1,\n\}$.\\[2mm]
Let $\Upsilon$ be a set of signed formulas. Then {\em
$h$ satisfies $\Upsilon$} if $h$ satisfies $\eta$ for every $\eta
\in \Upsilon$.

\begin{lem} Let $h: \mathfrak{Fm}'' \to \mathfrak{M}_{4m}$ be an
homomorphism, and let
$$
\begin{array}{l} \displaystyle \frac{\eta}
{\Upsilon_1 \ | \ \ldots \ | \ \Upsilon_n}
\end{array}
$$
be a rule of $\mathbb{T}$. If $h$ satisfies $\eta$ then $h$
satisfies $\Upsilon_i$ for some $1 \leq i \leq n$.
\end{lem}
\begin{dem}
A straightforward proof by cases.
\end{dem}

\

\begin{prop} \label{prop5.4} If \ $\not\models_{\tml}
\alpha$ \ then every completed tableau for $F(\alpha)$ is open.
\end{prop}
\begin{dem}
Assume that \ $\not\models_{\tml} \alpha$ \ and suppose
that there exists a completed closed tableau $\mathcal{T}$ for
$F(\alpha)$. Since $\not\models_{\tml} \alpha$, there is an
homomorphism $h$ such that $h(\alpha) \in \{\0, \n\}$. Thus, $h$
satisfies $F(\alpha)$. By the previous lemma, $h$ must satisfy
the set of signed formulas occurring in some branch $\theta$ of
$\mathcal{T}$. Since $\mathcal{T}$ is closed, the branch $\theta$
is closed, that is, the closure rule was used in $\theta$. But
it is an easy task to verify that no homomorphism can satisfy  simultaneously both premises of the
closure rule. This leads to a contradiction, and then every
completed tableau for $F(\alpha)$ must be open.
\end{dem}

\

\begin{cor}[Soundness of $\mathbb{T}$] \label{soundtabl} If \
$\vdash_{\mathbb{T}} \alpha$ then \ $\models_{\tml} \alpha$.
\end{cor}

\noindent In order to prove completeness, we need to state the
following result.

\begin{prop} \label{hom-table}
Let $\theta$ be an open branch of a completed  open tableau
$\mathcal{T}$, and let $\Upsilon$ be the set of signed formulas
occurring in $\theta$. Let $h$ be an homomorphism such that, for
every $\alpha \in Var$:\\

$$(\ddagger\ddagger) \left \{ \begin{tabular}{rl}
$h(\alpha) \in \{\1, \b\}$  & if \ $T(\alpha) \in \Upsilon$;\\
$h(\alpha) \in \{\1,\n\}$  & if \ $F(\neg\alpha) \in \Upsilon$;\\
$h(\alpha) \in \{\0,\n\}$  & if \ $F(\alpha) \in \Upsilon$;\\
$h(\alpha) \in \{\0,\b\}$  & if \ $T(\neg\alpha) \in \Upsilon$.\\
\end{tabular}\right.$$

\noindent In any other case $h(\alpha)$ is arbitrary, for $\alpha
\in Var$. Then $(\ddagger\ddagger)$ holds for any complex formula
$\alpha$.
\end{prop}
\begin{dem}
By induction on the degree of $\alpha$.\\
(i) $\alpha \in Var$. The result is clearly true.\\
(ii) $\alpha=\neg\beta$. If $T(\alpha) \in \Upsilon$ then
$T(\neg\beta) \in \Upsilon$ and so $h(\beta) \in \{\0,\b\}$, by
induction hypothesis. Thus, $h(\alpha) \in \{\1,\b\}$. If $T(\neg\alpha) \in \Upsilon$ then
$T(\neg\neg\beta) \in \Upsilon$ and so $T(\beta) \in \Upsilon$, since $\mathcal{T}$ is completed. Thus $h(\beta) \in \{\1,\b\}$, by
induction hypothesis and then $h(\alpha) \in \{\0,\b\}$. The other
cases are proved analogously.\\
(iii) $\alpha=\beta \succ\gamma$.\\
(iii.1) $T(\alpha) \in \Upsilon$. Since $\mathcal{T}$ is completed
then the rule for $T(\beta \succ\gamma)$ was eventually used,
splitting into three branches. One of them is a sub-branch of
$\theta$, thus one of the following cases hold:\\
(iii.1.1) $T(\gamma)\in \Upsilon$. Then $h(\gamma\in \{\1,\b\}$, by the inductive hypothesis. Then, $h(\alpha)\in \{\1,\b\}$, by the definition of $\succ$. \\
(iii.1.2) $T(\neg\beta), F(\gamma), T(\neg\gamma) \in \Upsilon$.
Then $h(\beta) \in \{\0,\b\}$ and $h(\gamma)=\0$, by induction
hypothesis, and so $h(\alpha) \in \{\1,\b\}$.\\
(iii.1.3) $F(\beta), F(\gamma), F(\neg\gamma) \in \Upsilon$. Then $h(\beta) \in\{\0, \b\}$,
by induction hypothesis, and so  $h(\alpha)=\1 \in \{\1,\b\}$.\\
The proof of the remaining sub cases for (iii) (namely: $F(\alpha) \in
\Upsilon$, $T(\neg\alpha) \in \Upsilon$ and $F(\neg\alpha) \in
\Upsilon$) are analogous.
\end{dem}

\

\begin{prop} \label{lem-comple}
Assume that there is a completed  open tableau for $F(\alpha)$.
Then \ $\not\models_{\tml} \alpha$.
\end{prop}
\begin{dem}
Consider, by hypothesis, an open branch $\theta$ of a completed
open tableau $\mathcal{T}$ for $F(\alpha)$, and let $\Upsilon$ be
the set of signed formulas occurring in $\theta$. Let $h$ be an
homomorphism defined as in Proposition~\ref{hom-table}. Then
$h(\alpha) \in \{\0,\n\}$, since $F(\alpha) \in \Upsilon$, and so  \
$\not\models_{\tml} \alpha$.
\end{dem}

\

\begin{theo}[Completeness of $\mathbb{T}$] \label{t-comple}
If \ $\models_{\tml} \alpha$ then \ $\vdash_{\mathbb{T}}
\alpha$.
\end{theo}
\begin{dem}
If \ $\models_{\tml} \alpha$ then, by
Proposition~\ref{lem-comple}, every completed tableau for
$F(\alpha)$ is closed, and so there exists (by
Proposition~\ref{completabl}) a closed tableau for $F(\alpha)$.
That is, \ $\vdash_{\mathbb{T}} \alpha$.
\end{dem}

\

\begin{cor} Let $\alpha$ be a formula. Then every completed tableau for
$F(\alpha)$ is open, or  every completed tableau for $F(\alpha)$
is closed.
\end{cor}
\begin{dem}
Suppose that there exists a completed open tableau $\mathcal{T}$
for $F(\alpha)$, as well as a completed closed tableau
$\mathcal{T}'$ for $F(\alpha)$. By Proposition~\ref{lem-comple},
\ $\not\models_{\tml} \alpha$. On the other hand, by
Corollary~\ref{soundtabl}, it follows that \ $\models_{\tml}
\alpha$, a contradiction.
\end{dem}

\

\begin{prop} \label{tabl-clos}
Suppose that $\vdash_{\mathbb{T}} \alpha$. Then every completed
tableau for $T(\neg\alpha)$ is closed.
\end{prop}
\begin{dem}
If $\vdash_{\mathbb{T}} \alpha$ then \ $\models_{\tml}
\alpha$, by Corollary~\ref{soundtabl}. Thus, $h(\alpha) \in
\{\1,\b\}$ for every homomorphism $h$. Suppose that there exists a
completed open tableau $\mathcal{T}$ for $T(\neg\alpha)$, and let
$\Upsilon$ be the set of signed formulas obtained from an open
branch $\theta$ of $\mathcal{T}$. Define an homomorphism $h$ as in
Proposition~\ref{hom-table}. Then $h(\alpha) \in \{\0,\b\}$, since
$T(\neg\alpha) \in \Upsilon$. But $h(\alpha) \in \{\1,\b\}$ and so
$h(\alpha)=\b$. Using Lemma~\ref{Lemhh'} there exists an
homomorphism $h'$ such that $h'(\alpha)=\n$, a contradiction.
Therefore every completed tableau for $T(\neg\alpha)$ is closed.
\end{dem}

\

\noindent The last result shows from the tableau perspective the
fact that, in $\tml$, the tautologies just get the
truth-value $\1$ by means of homomorphisms. Moreover,
Proposition~\ref{tabl-clos} will allow to prove, only by tableau
tools, the admissibility of the Rule of Necessitation (Nec). In order
to see this, suppose that \ $\models_{\tml} \alpha$, and
start a tableau in $\mathbb{T}$ for $F(\square\alpha)$. By
Proposition~\ref{der-rule2}, two branches are originated: one
with $F(\alpha)$ and the other with $T(\neg\alpha)$. By
Proposition~\ref{completabl}, both tableaux will eventually
terminate. Using Theorem~\ref{t-comple} and
Proposition~\ref{tabl-clos} both tableaux are closed and so the
original tableau for $F(\square\alpha)$ is closed. This shows that
\ $\models_{\tml} \square\alpha$, by
Corollary~\ref{soundtabl}.

\

\noindent It is worth noting that the tableau system
$\mathbb{T}$ allows to decide whether a given formula is valid or
not in $\tml$, and so it decides the validity in the
variety ${\bf TMA}$ of equations of the form $\alpha
\approx \1$. With respect to inferences of the form
$\alpha\vdash_{\tml} \beta$, they can be recovered in
$\mathbb{T}$ by means of tableaux for $F(\alpha \succ \beta)$.
Thus, $\mathbb{T}$ decides the validity in the variety ${\bf
TMA}$ of equations $\alpha \approx \beta$. Finally, as
it happens with the classical case (cf.~\cite{Smu}), the set of
signed formulas of an open branch of a completed open tableau
allows to find a model for that set of formulas: in particular,
it finds a counter-model for a non-valid formula. 

\section{A tableau system for $\tml$ in the original language}\label{s6}

In this short section, we use the results exposed in the previous section to present a tableau system for $\tml$ in the original language.

\begin{defi} \label{der-rule2} Let $\mathbb{T}^{\prime}$ be the tableau system defined by the following set of rules:\\

$$ \displaystyle \frac{T(\alpha \vee\beta)}{ T(\alpha) \, | \, T(\beta)} \hspace{1cm} \frac{T(\neg(\alpha \vee\beta))}{ T(\neg\alpha), T(\neg\beta)}$$

\

$$ \displaystyle \frac{F(\alpha \vee\beta)}{F(\alpha), F(\beta)} \hspace{1cm} \frac{F(\neg(\alpha \vee\beta))}{ F(\neg\alpha) \, | \, F(\neg \beta)}$$

\

$$ \displaystyle \frac{T(\alpha \wedge\beta)}{T(\alpha), T(\beta)}\hspace{1cm} 
\frac{T(\neg(\alpha \wedge\beta))}{T(\neg\alpha) \, | \,  T(\neg\beta)}$$

\

$$ \displaystyle \frac{F(\alpha \wedge\beta)}{F(\alpha) \, | \,  F(\beta)} \hspace{1cm} 
\frac{F(\neg(\alpha \wedge\beta))}{F(\neg\alpha), F(\neg\beta)}$$

\

$$\begin{array}{ccc}
\displaystyle \frac{T(\neg\neg\alpha)}{T(\alpha)} &
\mbox{\ \ \ \ } & \displaystyle
\frac{F(\neg\neg\alpha)}{F(\alpha)}
\end{array}$$

\

$$
\begin{array}{cccc}
\begin{array}{l} \displaystyle \frac{T(\square\alpha)}
{T(\alpha), F(\neg\alpha)}
\end{array}
&
\begin{array}{l} \displaystyle \frac{F(\square\alpha)}
{F(\alpha)\ | \ T(\neg\alpha)}
\end{array}
&
\begin{array}{l} \displaystyle \frac{T(\neg\square\alpha)}
{F(\square\alpha)}
\end{array}
&
\begin{array}{l} \displaystyle \frac{F(\neg\square\alpha)}
{T(\square\alpha)}
\end{array}

\end{array}
$$

\

\noindent {\bf Closure rule:}
$$
\displaystyle \frac{T(\alpha), F(\alpha)}{\star}$$

\
\end{defi}

\

The satisfaction of a signed formula $\eta$ by a homomorphism $h: \mathfrak{Fm}\to \mathfrak{M}_{4m}$ is defined as in Subsection \ref{s5.4}. Then, 

\begin{lem} Let $h: \mathfrak{Fm}\to \mathfrak{M}_{4m}$ be an
homomorphism, and let
$$
\begin{array}{l} \displaystyle \frac{\eta}
{\Upsilon_1 \ | \ \ldots \ | \ \Upsilon_n}
\end{array}
$$
be a rule of $\mathbb{T}^{\prime}$ for $n=0,1,2$. If $h$ satisfies $\eta$ then $h$
satisfies $\Upsilon_i$ for some $1 \leq i \leq n$.
\end{lem}
\begin{dem}
A straightforward proof by cases.
\end{dem}

\begin{prop} If \ $\not\models_{\tml}
\alpha$ \ then every completed tableau for $F(\alpha)$ is open.
\end{prop}
\begin{dem} Analogous to the proof of Proposition \ref{prop5.4}.
\end{dem}

\begin{cor}[Soundness of $\mathbb{T}^{\prime}$] If \
$\vdash_{\mathbb{T}} \alpha$ then \ $\models_{\tml} \alpha$.
\end{cor}

\noindent Next, we prove a version of Proposition \ref{hom-table} for this setting.

\begin{prop} \label{hom-table1}
Let $\theta$ be an open branch of a completed  open tableau
$\mathcal{T}$, and let $\Upsilon$ be the set of signed formulas
occurring in $\theta$. Let $h$ be an homomorphism such that, for
every $\alpha \in Var$:\\

$$(\ddagger\ddagger) \left \{ \begin{tabular}{rl}
$h(\alpha) \in \{\1, \b\}$  & if \ $T(\alpha) \in \Upsilon$;\\
$h(\alpha) \in \{\1,\n\}$  & if \ $F(\neg\alpha) \in \Upsilon$;\\
$h(\alpha) \in \{\0,\n\}$  & if \ $F(\alpha) \in \Upsilon$;\\
$h(\alpha) \in \{\0,\b\}$  & if \ $T(\neg\alpha) \in \Upsilon$.\\
\end{tabular}\right.$$

\noindent In any other case $h(\alpha)$ is arbitrary, for $\alpha
\in Var$. Then $(\ddagger\ddagger)$ holds for any complex formula
$\alpha$.
\end{prop}
\begin{dem}
By induction on the complexity of $\alpha$ (see Definition \ref{defi01}). We analyze just some cases, the others are analogous.\\
(i) $\alpha \in Var$. The result is clearly true.\\
(ii) $\alpha=\neg\beta$. If $T(\alpha) \in \Upsilon$ then
$T(\neg\beta) \in \Upsilon$ and so $h(\beta) \in \{\0,\b\}$, by
induction hypothesis. Thus, $h(\alpha) \in \{\1,\b\}$. If $T(\neg\alpha) \in \Upsilon$ then
$T(\neg\neg\beta) \in \Upsilon$ and so $T(\beta) \in \Upsilon$, since $\mathcal{T}$ is completed. Thus $h(\beta) \in \{\1,\b\}$, by
induction hypothesis and then $h(\alpha) \in \{\0,\b\}$. The other
cases are proved analogously.\\
(iii) $\alpha=\beta \wedge\gamma$.\\
Suppose that $T(\alpha) \in \Upsilon$. Since $\mathcal{T}$ is completed
then the rule for $T(\beta \wedge\gamma)$ was eventually used. Then $T(\beta), T(\gamma) \in \Upsilon$, and by the inductive hypothesis, $h(\beta), h(\gamma) \in \{\1,\b\}$ and therefore $h(\beta\wedge\gamma)\in  \{\1,\b\}$.\\
If $T(\neg\alpha) \in \Upsilon$, then the rule for $T(\neg(\beta \wedge\gamma))$ was eventually used splitting into two branches. One of them is a sub-branch of
$\theta$, thus one of the following conditions hold:\\
(iii.1) $T(\neg\beta)\in \Upsilon$. Then $h(\neg\beta)\in \{\1,\n\}$, by the inductive hypothesis. Then, $h(\neg\alpha)\in \{\1,\n\}$, by the definition of of the operations in $\tml$ and by the fact that $h(\neg \alpha)=h(\neg\beta\vee\neg\gamma)$. \\
(iii.2) $T(\neg\gamma)\in \Upsilon$. Analogous to  (iii.1).\\
The proof of the remaining sub cases, namely $F(\alpha) \in
\Upsilon$ and $F(\neg\alpha) \in \Upsilon$), are analogous. \\
(iv) $\alpha=\square \alpha$. \\
We just analyze the case where $F(\square \alpha)\in \Upsilon$. Since $\mathcal{T}$ is completed, the rule $T(\square \alpha)$ was used splitting into two branches. Then, one of the following conditions hold.\\
(iv.1) $F(\alpha) \in \Upsilon$. By the inductive hypothesis, $h(\alpha)\in\{\0, \n\}$ and then $h(\square\alpha)=0 \in \{\0, \n\}$.\\
(iv.2) $T(\neg\alpha)\in \Upsilon$.   By the inductive hypothesis, $h(\alpha)\in\{\0, \b\}$ and therefore $h(\square\alpha)=0\in \{\0, \n\}$. 
\end{dem}

\

\noindent Finally,

\begin{theo}[Completeness of $\mathbb{T}^{\prime}$] \label{t-compleN}
If \ $\models_{\tml} \alpha$ then \ $\vdash_{\mathbb{T}^{\prime}}
\alpha$.
\end{theo}
\begin{dem} Let $\alpha$ be a formula. Note that (as in Proposition \ref{lem-comple}) if there is a complete open tableau for $F(\alpha)$ then $\not\models_{\tml} \alpha$. So, if $\models_{\tml} \alpha$ then every tableau for $F(\alpha)$ is closed and therefore there exists a closed tableau for $F(\alpha)$. That is, $\vdash_{\mathbb{T}^{\prime}} \alpha$.
\end{dem}

\section{Concluding remarks}
In this paper, we continue with the study of $\tml$ under a proof-theoretic perspective. In first place, we show that the natural deduction system ${\bf ND}_{\cal TML}$ introduced \cite{MF} admits a normalization theorem. Later, taking advantage of the contrapositive implication for the tetravalent modal algebras introduced in \cite{AVF3}, we define a decidable tableau system for $\tml$. The original language of the logic of TMAs -- in
particular, the language of logic $\tml$-- does not have an implication connective as a primitive connective. However, using  the contrapositive implication for $\tml$ as a primitive connective and following a general
techniques introduced in~\cite{CCCM}, we define a sound an complete tableau system for this logic. Finally, inspired in this last system, we provide a sound and complete tableau system for $\tml$ in the original language. Finally, inspired in this last system, we provide a sound and complete tableau system for $\tml$ in the original language.
These two tableau systems constitute new (proof-theoretic) decision procedures for checking validity in the variety of tetravalent modal algebras, besides the four-valued truth-tables of $\tml$ and the one available in terms of the cut-free sequent system introduced in \cite{MF}.

An interesting task (for future works) would be to provide a natural deduction system for $\tml$ in terms of negation and implication ($\neg$, $\succ$). This would involve to find proper introduction/elimination rules for $\succ$ which does not seem easy an easy job taking into account the unusual properties of this implication.

Also, we propose to extend $\tml$ to first order languages. This would provide a suitable context for studying and developing its potential applications in Computer Science as envisaged by Antonio Monteiro forty-five years ago.


\begin{thebibliography}{99}

\bibitem{AnBel} Anderson, A. R.  and Belnap N. D. (with contributions by thirteen others), {\em Entailment: the logic of relevance and necessity}, volume II, (1992)  Princeton University Press.  

\bibitem{Avr}   Arieli, O. and Avron, A., The value of the four values. Artificial Intelligence
v. 102, n. 1 (1998), pp. 97--141.

\bibitem{Avron02} Avron, A., Ben-Naim, J. and  Konikowska, B., {\em Cut-free ordinary sequent calculi for logics having generalzed finite--valued semantics}. Logica Universalis, 1,  41--69, 2006.

\bibitem{Bel} Belnap, N., How computers should think. In: {\em Contemporary Aspects of Philosophy} (Editor: G. Ryle). Oriol Press, pp. 30--56, 1976.

\bibitem{Bez} B\'eziau, J.-Y., A new four-valued approach to modal logic.
{\em Logique et Analyse}, v. 54, n. 213 (2011).



\bibitem{ModLog} Blackburn, P., de Rijke, M. and Venema, Y., Modal Logic.
Cambridge University Press, 2001. ISBN 0-521-80200-8

\bibitem{CCCM}   Caleiro, C., Carnielli, W.A., Coniglio, M.
E. and Marcos, J. Two's company: The humbug of many logical
values. In: {\em Logica Universalis} (Editor J.-Y. B\'eziau).
Basel: Birkh\"auser, pp. 169-189, 2005.


\bibitem{CalMar} Caleiro, C. and Marcos, J., Classic-Like Analytic Tableaux for
Finite-Valued Logics. In: {\em
Logic, Language, Information And Computation}, Lecture Notes in
Computer Science vol. 5514, pp. 268-280. Eds.: H. Ono; M.
Kanazawa and R. de Queiroz. Springer, 2009.

\bibitem{CaFi} Cantú, L. M. and Figallo, M., Cut-free sequent-style systems for a logic associated to involutive Stone algebras. First online in {\em Journal of Logic and Computation} (2022).   https://doi.org/10.1093/logcom/exac061,

\bibitem{CCGGS} Carnielli, W.A., Coniglio, M.E., Gabbay, D., Gouveia, P. and Sernadas, C., Analysis and Synthesis of Logics. Volume 35 in the Applied Logic Series, Springer, 2008. 

\bibitem{WCMCJM} Carnielli, W.A., Coniglio, M.E. and Marcos, J., Logics of Formal
Inconsistency.
In: {\em Handbook of Philosophical Logic}, vol. 14, pp. 15-107.
Eds.: D. Gabbay; F. Guenthner. Springer, 2007.


\bibitem{Tax} Carnielli, W.A. and Marcos, J., A taxonomy of {{\bf C}}-systems.
In W.~A. Carnielli, M.~E. Coniglio, and I.~M.~L. D'Ottaviano,
editors, {\em Paraconsistency --- The logical way to the
inconsistent}, volume 228 of {\em Lecture Notes in Pure and
Applied Mathematics}, pp. 1--94.  Marcel Dekker, New York, 2002.

\bibitem{CF} M. Coniglio and M. Figallo, {\em Hilbert-style Presentations of Two Logics Associated  to Tetravalent Modal Algebras }, Studia Logica, nro. 3, vol. 102 (2014), 525--539.

\bibitem{daC} Da Costa, N.C.A., \emph{{I}nconsistent {F}ormal {S}ystems} (in
{P}ortuguese). {H}abilitation {T}hesis, 1963. Republished by
Editora UFPR, Curitiba, 1993.

\bibitem{daC2} Da Costa, N.C.A., Calculs propositionnel pour les syst\`emes
formels inconsistants. {\em Comptes Rendus de l'Acad\'emie de
Sciences de Paris}, s\'erie A, vol. 257(1963), 3790--3792.

\bibitem{MF} Figallo, M. {\em Cut-free Sequent Calculus and Natural Deduction for the Tetravalent Modal Logic}. Studia Logica, {\bf 109\/}, 1347--1373.(2021). https://doi.org/10.1007/s11225-021-09944-3


\bibitem{AVF3} Figallo, A.V. and Landini, P., On generalized I-algebras and 4-valued
modal algebras. {\em Reports on Mathematical Logic} 29 (1995),
3--18.

\bibitem{AVF2} Figallo, A.V. and Ziliani, A., Symmetric tetra-valued modal algebras.
{\em Notas Soc. Mat. Chile}  v. 10, n. 1 (1991), 133--141.

\bibitem{FR1} Font, J.M. and Rius, M., A four-valued modal logic arising from
Monteiro's last algebras. In {\em Proc. 20th Int. Symp.
Multiple-Valued Logic} (Charlotte, 1990), The IEEE Computer
Society Press, 85--92, 1991.

\bibitem{FR2} Font, J.M. and Rius, M., An abstract algebraic logic approach to
tetravalent modal logics. {\em  J. Symbolic Logic} v. 65, n. 2
(2000), 481--518.

\bibitem{GG} Gastaminza, M.L.  and Gastaminza, S., Characterization of a De Morgan lattice in terms of implication and negation. Proc. Japan Acad. v. 44, n. 7 (1968), 659--662.

\bibitem{Jans} Jansana, R., Propositional Consequence Relations and Algebraic Logic, {\em The Stanford Encyclopedia of Philosophy} (Spring 2011 Edition), E.N. Zalta (ed.).\\ 
{\small http://plato.stanford.edu/archives/spr2011/entries/consequence-algebraic/}


\bibitem{LeSc} Lemmon, E.J. and Scott, D., An Introduction to Modal Logic. In:
{\em The Lemmon notes}, K. Segerberg (editor). Volume 11 of
American Philosophical Quarterly Monograph series. Basil
Blackwell, Oxford, 1977.

\bibitem{L1} Loureiro, I., {\em \'Algebras Modais Tetravalentes}. PhD thesis,
Faculdade de Ci\^encias de Lisboa, 1983.

\bibitem{L2} Loureiro, I., Homomorphism kernels of a tetravalent modal algebra.
{\em Portugaliae Mathematica}, 39 (1980), 371--377.

\bibitem{MoRo} Montgomery, H. and R. Routley, Contingency and non-contingency bases
for normal modal logics. {\em Logique et analyse}, 9 (1966), 318--328.

\bibitem{OdWa01}  Odintsov, S. P. and Wansing H.. Modal logics with Belnapian truth values. {\em Journal of Applied
Non-Classical Logics} , 20(3):279–301, 2010.

\bibitem{OdWa02}  Odintsov, S. P. and Wansing H.. Disentangling FDE-based paraconsistent modal logics. {\em Studia
Logica}, 105(6):1221–1254, 2017.

\bibitem{Pri} Priest, G.. Many-valued modal logics: a simple approach. {\em The Review of Symbolic Logic}, 1(02),
2008.

\bibitem{RiJuJa} Rivieccio, U, Jung A., and Jansana, R.. Four-valued modal logic: Kripke semantics and duality.
{\em Journal of Logic and Computation}, 27(1):155–199, jun 2015.

\bibitem{Smu} Smullyan, R.M., First-order logic. Berlin: Springer-Verlag,
1968. Corrected republication by Dover Publications, New York,
1995. 

\bibitem{TS} Troelstra, A. S. and Schwichtenberg,  H.,   {\em  Basic Proof System}. Cambridge, UK: Cambridge University Press (1996).

\bibitem{AZ} Ziliani, A., {\em Algebras de De Morgan modales $4$-valuadas mon\'adicas}.
PhD thesis, Universidad Nacional del Sur (Bah\'{\i}a Blanca),
2001.


\end{thebibliography}
\end{document}